\documentclass{jfm}

\usepackage{natbib}

\usepackage{amsmath,amsfonts,amssymb}
\usepackage{latexsym}
\usepackage{epsfig,overpic}
\usepackage{subfigure}
\usepackage{color}
\usepackage{ulem}

\newcommand{\bn}{\mathbf n}

\newcommand{\bu}{\mathbf u}
\newcommand{\bv}{\mathbf v}

\newcommand{\bx}{\mathbf x}

\DeclareGraphicsExtensions{.pdf,.eps,.ps,.eps.gz,.ps.gz,.eps.Y}

\def\eps{\varepsilon}
\newcommand{\lbl}[1]{\label{#1}}

\newcommand{\ignore}[1]{}

\newcommand{\bm}{\boldsymbol}

\newcommand{\none}{{{n+1}}}

\newcommand{\wtilde}{\widetilde}

\newcommand{\Bn}{{\bm n}}

\newcommand{\Grad}[1]{\nabla #1}

\newcommand{\veps}{\varepsilon}

\begin{document}
\title[Sharp-interface limits of a phase-field model for moving contact lines]{Sharp-interface limits of a phase-field model with a generalized Navier slip boundary condition for moving contact lines}
\author[X. Xu, Y. Di and H. Yu]{Xianmin XU, Yana DI \and Haijun YU\thanks{Email address for correspondence: hyu@lsec.cc.ac.cn}}

\affiliation{
School of Mathematical Sciences,
	University of Chinese Academy of Sciences, Beijing
	100049, China\\ 
	NCMIS \& LSEC, Institute of
	Computational Mathematics and Scientific/Engineering
	Computing, Academy of Mathematics and Systems
	Science, Beijing 100190, China}
\date{11 October 2017}
\maketitle
\begin{abstract}
The sharp-interface limits of a phase-field model with a generalized Navier slip boundary condition for moving contact line problem
are studied  by asymptotic analysis and numerical simulations. The effects of the {mobility} number 
as well as a phenomenological relaxation parameter in the boundary condition are considered. 
In asymptotic analysis, we focus on the case that the {mobility} number is the same order of the Cahn number
and derive the sharp-interface limits for several setups of the boundary relaxation parameter.
It is shown that the sharp interface limit of the phase field model is the standard two-phase 
incompressible Navier-Stokes equations coupled with several different slip boundary conditions.
Numerical results are consistent with the analysis results and also illustrate 
the different convergence rates of the sharp-interface limits for different scalings of the two parameters.
\end{abstract}

\section{Introduction}
Moving contact lines are common in nature and our daily life, 
e.g. the motion of rain drops on window glass, coffee rings left by 
evaporation of coffee drops, wetting on lotus leaves, etc.
The moving contact line problem also has many applications in {some} industrial processes,
like painting, coating and oil recovery, etc. Therefore, the problem has been studied extensively. see recent
review papers by \citet{pismen2002mesoscopic,blake2006physics,bonn2009wetting,snoeijer2013moving} and the references therein.

Moving contact line is a  challenging problem in fluid dynamics.
The standard two-phase Navier-Stokes
equations with a no-slip boundary condition will lead to a non-physical non-integrable
stress~\citep{huh1971hydrodynamic,dussan1979spreading}. This is the so-called 
contact line paradox. There are many efforts to solve {this} paradox. A natural way is to 
relax the no-slip boundary condition.
Instead, one could use the Navier slip boundary condition~\citep{huh1971hydrodynamic,cox1986,zhou1990dynamics,haley1991effect, spelt2005level,ren2007boundary}. In some applications, an effective
slip condition can be induced by  numerical methods~\citep{renardy2001numerical,marmottant2004spray}. The other approaches
to cure the paradox
include: to assume a precursor thin film and a disjoint pressure~\citep{schwartz1998simulation,pismen2000disjoining,eggers2005contact};
to derive a new thermodynamics for surfaces~\citep{shikhmurzaev1993moving}; to treat the
moving contact line as  a thermally activated process~\citep{blake2006physics,blake2011dynamics,seveno2009dynamics}, to use a
diffuse interface model for moving contact lines~\citep{seppecher1996moving,gurtin1996two,Jacqmin2000,QianWangSheng2003,yue2011can}, etc.

The diffuse interface approach for moving contact lines {has} become popular  recent years~\citep{anderson1998diffuse,QianWangSheng2004,ding2007wetting,carlson2009modeling,ren2011derivation,sibley2013contact,sui2014numerical,shen2015efficient,fakhari2017diffuse}.
In this approach, the interface is a thin diffuse layer between different fluids represented by a phase field function.
Intermolecular diffusion, caused by
the non-equilibrium of the chemical potential, occurs in the thin layer. The chemical 
diffusion can cause the motion of the contact line, even without using a slip boundary condition
on the solid boundary~\citep{Jacqmin2000,chen2000interface,briant2004lattice,Yue2011,Sibley2013}. 
On the other hand,  it is possible to 
combine the diffuse interface model with some slip boundary condition.
\citet*{QianWangSheng2003} proposed a
phase-field model with a generalized Navier slip boundary 
condition {(GNBC)}. The 
model takes account of the effect of the uncompensated Young stress,
which is  important to understand the difference  of
the dynamic contact angle and the Young's angle in molecular scale~\citep{QianWangSheng2003,ren2007boundary}.
Theoretically, the model can be derived from the Onsager variational
principle~\citep{QianWangSheng2006}. 
Numerical simulations {using this} model fit remarkably well with the molecular dynamics simulations~\citep{QianWangSheng2003}
and physical experiments~\citep{guo2013direct}. 
The model has also been used in problems with chemically patterned boundaries~\citep{WangQianSheng2008}, dynamic
wetting problems~\citep{carlson2009modeling,yamamoto2014modeling}, etc. {Several numerical} methods for the model have been 
developed~\citep{gao2012gradient,bao2012finite,gao2014efficient,shen2015efficient,aland2015efficient}.

Phase field models are convenient for  numerical calculations
~\citep{yue2004diffuse,teigen2011diffuse,sui2014numerical}.
One does not need to track the interface explicitly
as in using a sharp interface model. The 
phase-field function, which usually described by a Cahn-Hilliard~\citep{cahn1958free} equation 
or an Allen-Cahn equation~\citep{allen1979microscopic}, can capture the interface implicitly
and automatically. 
This makes  computations and 
analysis for the phase field model much easier than other approaches.
However, there are also some restrictions to use a diffuse interface 
model in real simulations. A key issue is that the thickness
of the diffuse interface can not be chosen as small as the physical size~\citep{khatavkar2006scaling},
due to the restriction of the computational resources. One
often choose a much larger (than physical values) interface thickness parameter(or 
a dimensionless Cahn number) in simulations.
But, only when phase field
 model approximates a sharp-interface limit 
correctly, the numerical simulations by this model 
with relatively large Cahn number can be trustful and
compared with experiments quantitatively.
Therefore, it is very important to study the sharp-interface limit of a phase field 
model~\citep{caginalp1998convergence,chen2014analysis}. 

The sharp-interface limits of diffuse interface models for two-phase flow without moving
contact lines has been studied a lot, both theoretically
and numerically~\citep{lowengrub1998quasi,jacqmin1999calculation,khatavkar2006scaling,huang2009mobility,magaletti2013sharp,sibley2013unifying}. In comparison, there {are} much
less studies for the sharp-interface 
limit of the phase field models for moving contact lines~\citep{Yue2010,kusumaatmaja2016moving}.
One important progress is made by \citet*{Yue2010}. They studied
the sharp interface limit of a phase field model with a no-slip boundary condition and
found a surprising result that only when the mobility parameter (denoted as  $\mathsf{L_d}$) is of order $O(1)$, 
the phase field model has a sharp-interface 
limit as the Cahn number (denoted as $\eps$) goes to zero. Notice that the usual choice of the mobility parameter
is of order $O(\eps^\beta)$, $1 \leq \beta\leq 3$ for problems without moving contact lines~\citep{magaletti2013sharp}.  For the phase field model
with the generalized Navier slip boundary condition, the only study for its sharp interface
limit is done by \citet*{wang2007sharp}. They also assumed  the mobility parameter is of order $O(1)$.
Their asymptotic analysis shows that the sharp-interface limit of the model 
is  a Hele-Shaw flow coupled with a standard Navier-slip boundary condition.
So far, it is not clear what is the sharp-interface limit of a phase field model
for moving contact line problem under the standard choice for the mobility parameter, e.g. $\mathsf{L_d}=O(\eps)$.
This is the motivation of our study.

We study the sharp-interface limit of
the phase field model with the GNBC by asymptotic analysis and 
numerical simulations.
In asymptotic analysis, we assume that the  mobility number $\mathsf{L_d}$ is of order $O(\varepsilon)$
and consider several typical scalings of phenomenological {boundary} relaxation parameter $\mathsf{V_s}$ in the GNBC model.
We show that the sharp-interface limit  is a standard  two-phase
Navier-Stokes equations coupled with different slip boundary conditions for  different choices of $\mathsf{V_s}$.
When $\mathsf{V_s}=O(\eps^{\beta})$ with $\beta=0,-1$, 
we obtain a sharp-interface version of  the  GNBC. In the case $\mathsf{V_s}=O(1)$,
the velocity of the contact line is equal to the fluid velocity, while in the case $\mathsf{V_s}=O(\eps^{-1})$,
the velocity of the contact line is different from the fluid velocity due to the contribution of the chemical diffusion
on the boundary.
When $\mathsf{V_s}=O(\eps^{-2})$, we obtain the standard Navier slip boundary 
condition together {with} the condition that the dynamic contact angle is equal to
the static contact angle.
Numerical experiments for a  Couette flow show the different sharp-interface 
limits for the various choice of $\mathsf{L_d}$ and $\mathsf{V_s}$. 
Furthermore, numerical results also {reveal} the different convergence 
rates for different choices of the two parameters.
For very large relaxation parameter $\mathsf{V_s}=O(\eps^{-3})$,
the numerical results are very similar to the results by~\citet*{Yue2010}.

The structure of the paper is as follows. In Section 2, we introduce
 the phase field model with the GNBC and its non-dimensionalization. 
 In Section 3, the sharp-interface limits
of the phase field model with the GNBC are obtained for various
choice of $\mathsf{V_s}$ by asymptotic analysis. In Section 4,
we show the numerical experiments for a
Couette flow by a recent developed second order scheme.
Finally, some concluding remarks are given in Section 5.

%
\section{The phase field model with generalized Navier slip boundary condition}
A Cahn-Hilliard-Navier-Stokes (CHNS) system with the generalized Navier  boundary condition (GNBC) is proposed by~\citet*{QianWangSheng2003} to describe a two-phase flow with  moving contact lines. The CHNS system reads
\begin{equation}\label{e:CHNS}
\left\{
  \begin{array}{ll}
    \frac{\partial \phi}{\partial t}+\mathbf v \cdot\nabla\phi = M\Delta \mu, \qquad \qquad\mu=-K\Delta\phi-r(\phi-\phi^3), & \hbox{} \\
    \rho[\frac{\partial \mathbf v}{\partial t}+(\mathbf v\cdot\nabla)\mathbf v]=\mathbf F-\nabla p+\eta\Delta\mathbf{v}+\mu\nabla\phi, \quad  \nabla \cdot \mathbf v=0. & \hbox{}
  \end{array}
\right.
\end{equation}
The first equation is the Cahn-Hilliard equation. Here
$\phi$ is the  phase field function,  and $\mu$ is the chemical potential.
The  thickness of the diffuse interface is $\xi=\sqrt{K/r}$ and the fluid-fluid interface tension 
is given by $\gamma=2\sqrt{2}r\xi/3$. 
$M$ is a phenomenological mobility coefficient.
The second equation in \eqref{e:CHNS} is the incompressible Navier-Stokes equation for two-phase flow.
Here 
$\mu\nabla\phi$ describes the capillary force exerted to the fluids by the interface.
For simplicity, we assume that the two fluids have equal  density $\rho$ and viscosity $\eta$.

The  generalized Navier boundary condition  on the solid boundary is
\begin{eqnarray}\label{e:GNBC}
& \beta (v_{\tau}-v_w)=-\eta\partial_{n} v_{\tau}+L(\phi)\partial_{\tau}\phi,
& \quad v_n=0,\\
& L(\phi)=K\partial_n \phi+\frac{\partial \gamma_{wf}(\phi)}{\partial \phi},
& \gamma_{wf}(\phi)=-\frac{\gamma}{4}\cos\theta_s(3\phi-\phi^3).
\end{eqnarray}

Here $v_n$ and $v_{\tau}$ are respectively the normal fluid velocity and the tangential fluid velocity on the solid boundary. $v_w$ is the velocity of the boundary itself. We assume the wall only moves in a tangential direction.
{$\beta$ is a slip coefficient and the slip length is given as $l_s=\eta/\beta$.}
$\gamma_{wf}(\phi)$
 is the solid-fluid interfacial energy density (up to a constant) and $\theta_s$ is the
static contact angle. $L(\phi)\partial_{\tau} \phi$ represents the uncompensated Young stress.

In addition, the boundary conditions for the phase field $\phi$ and the chemical potential $\mu$ are given by
\begin{equation}\label{e:relaxationBC}
\begin{array}{l}
\frac{\partial \phi}{\partial t}+v_{\tau}\partial_{\tau}\phi=-\Gamma L(\phi), \\
 \partial_n\mu=0,
\end{array}
\end{equation}
with $\Gamma$ being a  positive phenomenological relaxation parameter. 

\begin{figure}
\vspace{-0.2cm}
\centering
   \includegraphics[width=4.in]{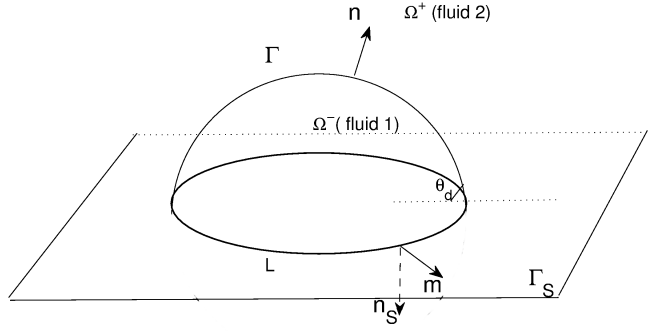}
\vspace{-1.5cm}
  \caption{A liquid drop on a planar solid surface $\Gamma_S$ with a contact line $L$.}
  \label{fig:1}
\end{figure}

To study the behavior of the CHNS system with the GNBC condition, it is useful to
  nondimensionalize the system.
Suppose the typical length scale in the two-phase flow system is given by $l$ and 
the characteristic velocity is $v^*$.
We then scale the velocity by $v^*$, the length by $l$,
the time by $l/v^*$, body force(density) $\mathbf F$ by $\eta v^*/l^2$ and the pressure by $  \eta v^*/l$.
With  six dimensionless parameters:
\begin{eqnarray}\label{e:dimensionless}
& \mathsf{L_d}=\frac{3M\gamma}{2\sqrt{2}v^*l^2} \hbox{ (the mobility number)}, & \mathsf{R_e}=\frac{\rho v^* l}{ \eta} \hbox{ (the Reynold number)}, \nonumber \\
&  \mathsf{B}=\frac{3\gamma}{2\sqrt{2}\eta v^*}\hbox{ (inverse of the Capillary number)}, &
 \mathsf{V_s}=\frac{3\gamma\Gamma l}{2\sqrt{2}v^*}\hbox{ (a relaxation parameter)}, \nonumber \\
&\mathsf{l_s}=\frac{ l_s}{ l}\hbox{ (the dimensionless slip length)}, & \eps =\frac{\xi}{l}\hbox{ (the Cahn number)},\nonumber
\end{eqnarray}
we have the following dimensionless Cahn-Hilliard-Navier-Stokes {system}
\begin{align}&\left\{
                  \begin{array}{ll}
\frac{\partial \phi}{\partial t}+\mathbf{v}\cdot\nabla\phi= \mathsf{L_d}\Delta\mu,\qquad
\mu=-\eps {\Delta}{\phi}-{\phi}/\eps+{\phi}^3/\eps,& \\
                     \mathsf{R_e}\Big[\frac{\partial\mathbf{v}}{\partial t}+(\mathbf{v}\cdot\nabla)\mathbf{v}\Big]=
\mathbf{F}  -\nabla{p} +{\Delta}{ \mathbf{v}} +  \mathsf {B} \mu {\nabla}{\phi}, & \hbox{} \\
                    {\nabla}\cdot{\mathbf{v}}=0, & \hbox{}
                  \end{array}
                \right. 
\lbl{e:2.1n}
\end{align}
with the boundary conditions
\begin{equation}
\left\{
  \begin{array}{ll}
     \frac{\partial {\phi}}{\partial {t}} + {{v}}_{\tau}{\partial_{\tau}}{\phi}=-\mathsf{V_s} \mathcal{L}(\phi), & \hbox{} \\
    \mathsf{l_s}^{-1} ({v}_{\tau} - v_w)=-  \partial_{n}{v}_{\tau}+ \mathsf{B}\mathcal{L}(\phi)\partial_{\tau}{\phi}, & \hbox{}\\
\partial_n\mu=0 , \ \ \ \   {v}_{n}=0, & \hbox{}
  \end{array}
\right.
\lbl{e:2.2n}
\end{equation}
where $\mathcal{L}(\phi)=\eps \partial_n{\phi} +\frac{\partial \gamma_{wf}(\phi)}{\partial \phi}$ and
$\gamma_{wf}(\phi)=-\frac{\sqrt{2}}{6}\cos\theta_s(3\phi-\phi^3)$ being the wall-fluid interfacial  energy density function. 
We {now}
 clarify some notations in the boundary condition. Suppose the unit outward normal vector on
the solid boundary is given by $\mathbf{n}_S$ (see Figure~\ref{fig:1}). Then, we have $v_n=\mathbf{v}\cdot\mathbf{n}_S$,
$v_{\tau}=\mathbf{v}-v_n\mathbf{n}_S$, $\partial_n=\mathbf{n}_S\cdot\nabla$ and $\partial_{\tau}=\nabla-\mathbf{n}_S(\mathbf{n}_S\cdot\nabla)$.

\section{The asymptotic analysis}
We  do asymptotic analysis for the Cahn-Hilliard-Navier-Stokes system \eqref{e:2.1n}-\eqref{e:2.2n}.
Here we assume the mobility number satisfies $\mathrm{L_d}=O(\eps)$.
 We show that such a choice of mobility will also lead to the standard 
two-phase Navier-Stokes equation inside the domain. 
Furthermore, this {assumption} also makes it possible
to derive proper boundary conditions for the sharp-interface limit of the diffuse-interface model.
We show that  different setups of  the relaxation parameter 
$\mathsf{V_s}=O(\eps^{\beta})$, $\beta=0,-1,-2$ will lead to different boundary conditions.

To make the presentation in this section clear, 
we  use $\phi_{\eps}$, $\mathbf{v}_{\eps}$
and $p_{\eps}$ instead of $\phi$, $\mathbf{v}$
and $p$ in the system \eqref{e:2.1n}-\eqref{e:2.2n}, to show explicitly 
that these functions depend on $\eps$.
We suppose that the system is located in a domain $\Omega$ with solid boundary $\Gamma_S{\subseteq} \partial\Omega$ (as shown in Figure~\ref{fig:1}).
Suppose that the two-phase interface is given by the zero level-set of the
 phase field function $\phi_\eps$
 \begin{equation}\label{e:interface}
 \Gamma:=\{x\in\Omega\ |\ \phi_\eps(x)=0\}.
 \end{equation}
{We denote by}  $\Omega^{-}=\{x\in\Omega\,|\, \phi_{\eps}<0\}$ the domain occupied by fluid 1  and 
$\Omega^{+}=\{x\in\Omega\,|\,\phi_{\eps}>0\}$ the domain occupied by fluid 2.

\subsection{The bulk equations}

We first do asymptotic analysis for the Cahn-Hilliard-Navier-Stokes system far from
the boundary. The analysis is the same as that for two-phase flow without 
contact lines. We will state the key steps of the analysis and illustrate the main results.
In the next subsection,
the bulk analysis here will be combined with the analysis near the boundary
to derive the sharp-interface limits of the GNBC.

Let  $\mathrm{L_d}=\eps l_d$. Consider the CHNS system far from the boundary
of $\Omega$.
We first do  outer expansions far from the interface $\Gamma$,
then we consider  inner expansions  near $\Gamma$. Combining them together,
we will obtain the sharp-interface limit of the  CHNS {system} in $\Omega$.

{\it Outer expansions. } 
Far from the two-phase interface $\Gamma$, we use the following ansatz,
\begin{equation}\lbl{e:outexpan}
\begin{array}{l}
\mathbf{v}_{\eps}^{\pm}=\mathbf{v}_0^{\pm}+\eps\mathbf{v}_1^{\pm}+\eps^2 \mathbf v_2^{\pm}+\cdots,\\
\phi_{\eps}^{\pm}=\phi_0^{\pm}+\eps\phi_1^{\pm}+\eps^2 \phi_2^{\pm}+\cdots,\\
p_{\eps}^{\pm}=p_0^{\pm}+\eps p_1^{\pm}+\eps^2 p_2^{\pm}+\cdots.
\end{array}
\end{equation}
Here $f^{\pm}$ denotes the restriction of a function $f$ in $\Omega^+$ and $\Omega^-$ respectively.
For $\mu_\eps$, we easily have
\begin{equation*}
\mu_\eps^{\pm}=\eps^{-1}\mu_0^{\pm}+\mu_1^{\pm}+\eps\mu_2^{\pm}+\cdots,
\end{equation*}
where
\begin{equation}\lbl{e:outer0}
\mu_0^{\pm}=-\phi_0^{\pm}+(\phi_0^{\pm})^3.
\end{equation}

We substitute the above expansions to the CHNS system \eqref{e:2.1n}.
The leading order of the first equation of \eqref{e:2.1n} gives
\begin{equation}\lbl{e:outer1a}
\frac{\partial \phi^{\pm}_0}{\partial t}+ \bv_0\cdot\nabla\phi_0^{\pm}=l_d \Delta \mu_0^{\pm}. 
\end{equation} 
The leading order of the second equation of \eqref{e:2.1n} gives
\begin{equation}\lbl{e:outer1b}
\mu_0^{\pm}\nabla\phi_0^{\pm}=0.
\end{equation}
More precisely, we have 
\begin{equation*}
(-\phi_0^{\pm}+(\phi_0^{\pm})^3)\nabla\phi_0^{\pm}=\nabla(\frac{(1-(\phi_0^\pm)^2)^2}{4})=0.
\end{equation*}
This implies that
\begin{equation}\lbl{e:phi0pm}
\phi_0^\pm=c_{\pm}\qquad \qquad \text{in } \Omega^{\pm},
\end{equation}
where $c_{\pm}$ are two constants such that $c_+>0$ and $c_-<0$.
For the third equation of \eqref{e:2.1n},  in the leading order, we have
\begin{align}\lbl{e:outer1c}
                    {\nabla}\cdot{\mathbf{v}_0^{\pm}}=0.
\end{align}
By direct calculations, we also have the next order of the second equation of \eqref{e:2.1n} as
\begin{equation}\lbl{e:outer2b}
 \mathsf{R_e}\Big[\frac{\partial\mathbf{v}_0^{\pm}}{\partial t}+(\mathbf{v}_0^{\pm}\cdot\nabla)\mathbf{v}_0^{\pm}\Big]=
\mathbf{F}  -\nabla{p}_0^{\pm} +{\Delta}{ \mathbf{v}_0^{\pm}} +\mu_0^{\pm}\nabla\phi_1^{\pm}.
\end{equation}
Here we have used the fact that $\phi_0^{\pm}$ are constants in $\Omega^{\pm}$.

{\it Inner expansions.} The outer expansion in $\Omega^{+}$ and $\Omega^-$
are connected by the transition layer near the interface $\Gamma$. We will
consider the so-called inner expansions near $\Gamma$. 
{Let} $d(x,t)$ be signed distance  to $\Gamma$, which is well-defined near
 the interface. Then
the unit normal of the interface pointing to $\Omega^+$ is given by $\mathbf{n}=\nabla d$.
We introduce a 
new rescaled  variable 
$$\xi=\frac{d(x)}{\eps}.$$
For any function $f(x)$ (e.g. $f=\mathbf{v}_\eps ,p_\eps, \phi_\eps$), we can rewrite it as
\begin{equation}
f(x)=\tilde{f}(x,\xi).
\end{equation}
Then we have
\begin{equation}\lbl{e:changevarible1}
\begin{array}{l}
\nabla f=\nabla\tilde{f}+\eps^{-1}\partial_\xi\tilde{f}\mathbf n,\\
\Delta f=\Delta \tilde{f}+\eps^{-1}\partial_\xi\tilde{f}\kappa+2\eps^{-1}(\bn\cdot\nabla)\partial_{\xi}\tilde{f}+\eps^{-2}\partial_{\xi\xi}\tilde{f},\\
\partial_t f=\partial_t\tilde{f}+\eps^{-1} \partial_\xi\tilde{f} \partial_t d_\eps.
\end{array}
\end{equation}
Here we use the fact that $\nabla\cdot\mathbf n=\kappa$,  the mean curvature of the 
interface. $\kappa(x)$ for $x\in\Gamma(t)$  is positive (resp. negative) if the domain $\Omega_-$ is convex (resp. concave) near $x$. 

In the inner region,  we  assume that
\begin{equation}\lbl{e:innerExpan}
\begin{array}{l}
\tilde{\mathbf{v}}_{\eps}=\tilde{\mathbf{v}}_0+\eps\tilde{\mathbf{v}}_1+\eps^2 \tilde{\mathbf v}_2+\cdots,\\
\tilde\phi_{\eps}=\tilde\phi_0+\eps\tilde\phi_1+\eps^2 \tilde\phi_2+\cdots,\\
\tilde p_{\eps}=\tilde p_0+\eps\tilde p_1+\eps^2\tilde p_2+\cdots.
\end{array}
\end{equation}
A direct expansion for the chemical potential $\tilde{\mu}_{\eps}$ gives
\begin{equation*}
\tilde\mu_\eps=\eps^{-1}\tilde\mu_0+\tilde\mu_1+\eps\tilde\mu_2+\cdots,
\end{equation*}
with
\begin{align}
&\tilde\mu_0=-\partial_{\xi\xi}\tilde{\phi}_0-\tilde\phi_0+\tilde\phi_0^3, \lbl{e:innerMu0}\\
&\tilde{\mu}_1=-\partial_{\xi\xi}\tilde{\phi}_1-\partial_\xi\tilde{\phi}_0\kappa+2(\bn\cdot\nabla)\partial_{\xi}\tilde{\phi}_0-\tilde\phi_1+3\tilde{\phi}_0^2\tilde{\phi}_1.
\lbl{e:innerMu1}
\end{align}

We substitute the above expansions into the system \eqref{e:2.1n}. Using the fact that $\mathrm{L_d}=\eps l_d$, 
in the leading order,
we have
\begin{equation}\lbl{e:inner1}
\left\{
\begin{array}{l}
\partial_{\xi\xi}\tilde{\mu}_0=0, \\
 \partial_{\xi\xi}\tilde{\mathbf{v}}_0+{\mathsf{B}}\tilde{\mu}_0\partial_\xi\tilde{\phi}_0\mathbf{n}=0,\\
 \mathbf{n}\cdot\partial_\xi\tilde{\mathbf{v}}_0=0.
\end{array}
\right.
\end{equation}
The next order  is 
\begin{equation}\lbl{e:inner2}
\left\{
\begin{array}{l} \partial_t d \partial_{\xi}\tilde{\phi}_0+ 
 \tilde{ \mathbf{v} }_0\cdot \mathbf{n}\partial_{\xi}\tilde{\phi}_0 =l_d(\partial_{\xi\xi}\tilde{\mu}_1+\kappa \partial_{\xi}\tilde{\mu}_0+2(\bn\cdot\nabla)\partial_{\xi}\tilde{\mu}_0 ),\\
 \tilde{\mathbf{v}}_0\cdot\mathbf{n}\partial_\xi\tilde{\mathbf{v}}_0=-\partial_\xi\tilde{p}_0\mathbf{n}+ \partial_{\xi\xi}\tilde{\mathbf{v}}_1+ \partial_\xi\mathbf{\tilde{v}}_0\kappa +2(\bn\cdot\nabla)\partial_{\xi}\tilde{\bv}_0+{\mathsf{B}}(\tilde \mu_1\partial_\xi\tilde{\phi}_0\mathbf{n}+\tilde \mu_0\partial_\xi\tilde{\phi}_1\mathbf{n}+\tilde \mu_0\nabla\tilde{\phi}_0),\\
 \mathbf{n}\cdot\partial_\xi\tilde{\mathbf{v}}_1+\nabla\cdot\tilde{\mathbf{v}}_0=0.
\end{array}
\right.
\end{equation}

{\it Matching conditions.} We need the following matching conditions for inner and outer expansions.
\begin{align}
&\lim_{\xi\rightarrow\pm\infty} \tilde{f}_i(x,\xi)= f_i^{\pm}(x), \label{e:match1}\\
&\lim_{\xi\rightarrow\pm\infty}  (\nabla_x \tilde{f}_i(x,\xi)+\partial_\xi\tilde f_{i+1}(x,\xi)\mathbf{n})=\nabla f_i^{\pm}(x).\label{e:match2}
\end{align}

In the following, we will derive the sharp-interface limit of the CHNS system \eqref{e:2.1n}
by the above inner and outer expansions.
From the first  equation of \eqref{e:inner1}, we know that $\tilde{\mu}_0$ is a linear function of $\xi$, which can be written as
$
\tilde{\mu}_0(\xi)=c_1\xi+c_0,
$
where $c_0$ and $c_1$ are independent of $\xi$. Since $\lim_{\xi\rightarrow\pm\infty}\tilde{\mu}_0=\mu^{\pm}$ is bounded,
we have $c_1=0$. Therefore 
\begin{equation}\lbl{e:TildeMu0}
\tilde{\mu}_0=c_0.
\end{equation}
Then the second  equation of \eqref{e:inner1} is reduced to
$$\partial_{\xi\xi}\tilde{\mathbf{v}}_0+{\mathsf{B}}c_0\partial_\xi\tilde{\phi}_0\mathbf{n}=0.$$
We integrate the equation with respect to $\xi$ in $(-\infty,\infty)$ and obtain
$$
\partial_{\xi}\tilde{\mathbf{v}}_0|_{-\infty}^{\infty}+{\mathsf{B}}c_0\tilde{\phi}_0|_{-\infty}^{\infty}\mathbf{n}=0.
$$
The inner product of the equation with $\mathbf{n}$ gives
$$
\mathbf{n}\cdot\partial_{\xi}\tilde{\mathbf{v}}_0|_{-\infty}^{\infty}+{\mathsf{B}}c_0\tilde{\phi}_0|_{-\infty}^{\infty}=0.
$$
By the third equation of \eqref{e:inner1}, we obtain that
$$ 
c_0\tilde{\phi}_0|_{-\infty}^{\infty}=0.
$$
By the matching condition, we have
$$ 
c_0(\phi_0^+-\phi_0^-)=0.
$$
Notice that $c_+=\phi_0^+>0>\phi_0^-=c_-$, we immediately have $c_0=0$, or equivalently 
\begin{equation}\lbl{e:mu0}
\tilde{\mu}_0=0.
\end{equation}
By the equation~\eqref{e:innerMu0}, we have 
\begin{equation}\lbl{e:Equprofile}
-\partial_{\xi\xi}\tilde{\phi}_0-\tilde\phi_0+\tilde\phi_0^3 =0.
\end{equation}
The solvability condition for this equation~\citep{pego1989}  leads to
\begin{equation}\lbl{e:solvability}
\lim_{\xi\rightarrow \pm\infty} \tilde{\phi}_0=\pm 1.
\end{equation}
And the solution of \eqref{e:Equprofile} is 
\begin{equation}\lbl{e:profile}
\tilde{\phi}_0=\tanh(\xi/\sqrt{2}).
\end{equation}
This is the profile of the $\tilde{\phi}_0$ in the diffuse-interface layer.
By the matching condition $\lim_{\xi\rightarrow \pm\infty} \tilde{\phi}_0=\phi_0^{\pm}$ and  \eqref{e:phi0pm}, we have
\begin{equation}
\phi_0^{\pm}(x)=c_{\pm}=\pm 1, \qquad\hbox{ in } \Omega^{\pm}.
\end{equation} 
This will lead to $\mu_0^{\pm}=0$.
Therefore the  equation~\eqref{e:outer2b} is reduced to
\begin{equation}
\mathsf{R_e}\Big(\frac{\partial\mathbf{v}_0^{\pm}}{\partial t}+(\mathbf{v}_0^{\pm}\cdot\nabla)\mathbf{v}_0^{\pm}\Big)=
\mathbf{F}  -\nabla{p}_0^{\pm} +{\Delta}{ \mathbf{v}_0^{\pm}}.
\end{equation}
Together with \eqref{e:outer1c}, this is the standard incompressible Navier-Stokes equation  in $\Omega^{\pm}$.

We then derive the jump conditions on the interface $\Gamma$. 
Noticing \eqref{e:mu0}, the second equation of \eqref{e:inner1} is reduced to
\begin{equation*}
\partial_{\xi\xi} \tilde{\mathbf{v}}_0=0.
\end{equation*}
{By} similar argument as for $\tilde{\mu}_0$ in \eqref{e:TildeMu0}, we know that
 $\tilde{\mathbf{v}}_0$ is independent of $\xi$, or 
 \begin{equation}\lbl{e:v0}
 \tilde{\mathbf{v}}_0(x,\xi)=\tilde{\mathbf{v}}_0(x).
 \end{equation} 
 By the matching condition 
 for $\tilde{\mathbf{v}}_0$, we obtain 
 \begin{equation}\lbl{e:jumpVel}
 [\mathbf{v}_0]=0,
 \end{equation}
where $[f]=f^{+}-f^{-}$ denotes the jump of a 
function $f$ across the interface $\Gamma$. The equation \eqref{e:jumpVel} implies that
the  velocity  is continuous across $\Gamma$.

Similarly, the first equation of \eqref{e:inner2} is reduced to  
\begin{equation*}
\partial_t d \partial_{\xi}\tilde{\phi}_0+ 
 \tilde{ \mathbf{v} }_0\cdot \mathbf{n}\partial_{\xi}\tilde{\phi}_0 =l_d \partial_{\xi\xi}\tilde{\mu}_1.
\end{equation*}
Integrate the equation in $(-\infty, \infty)$ and use the matching condition for $\tilde{\phi}_0$ and $\tilde{\mu}_1$, 
then we obtain
\begin{equation*}
\partial_t d+ \tilde{ \mathbf{v} }_0\cdot \mathbf{n}=0.
\end{equation*}
This implies that the normal velocity of the interface $\Gamma$ is 
\begin{equation}\lbl{e:normalVel}
V_n=  \mathbf{v}_0\cdot \mathbf{n}.
\end{equation}

We then show the jump condition for the viscous stress.
Noticing \eqref{e:mu0} and \eqref{e:v0},  the second equation of \eqref{e:inner2} can be reduced to
\begin{equation}\lbl{e:temp0}
-\partial_\xi\tilde{p}_0\mathbf{n}+ \partial_{\xi\xi}\tilde{\mathbf{v}}_1 +{\mathsf{B}}\mu_1\partial_\xi\tilde{\phi}_0\mathbf{n}=0.
\end{equation}
By the equation \eqref{e:innerMu1}, we have 
\begin{align}\lbl{e:temp1}
\int_{-\infty}^{\infty}\tilde{\mu}_1\partial_\xi\tilde{\phi}_0 d\xi=-\int_{-\infty}^{\infty}(\partial_\xi\tilde{\phi}_0)^2 d\xi \kappa+\int_{-\infty}^{\infty}(-\partial_{\xi\xi}\tilde{\phi}_1-\tilde\phi_1+3\tilde{\phi}_0^2\tilde{\phi}_1)\partial_\xi\tilde{\phi}_0 d\xi.
\end{align} 
Integration by parts for the second term in the right hand side of the equation leads to
 \begin{align*}
 \int_{-\infty}^{\infty}(-\partial_{\xi\xi}\tilde{\phi}_1-\tilde\phi_1+3\tilde{\phi}_0^2\tilde{\phi}_1)\partial_\xi\tilde{\phi}_0 d\xi=
 \int_{-\infty}^{\infty}(\partial_{\xi\xi}\tilde{\phi}_0+\tilde\phi_0-\tilde{\phi}_0^3)\partial_\xi\tilde{\phi}_1 d\xi=0.
 \end{align*}
 Here we use the equation~\eqref{e:Equprofile}. Then \eqref{e:temp1} is reduced to
\begin{align*}
\int_{-\infty}^{\infty}\tilde{\mu}_1\partial_\xi\tilde{\phi}_0 d\xi=-\sigma \kappa,
\end{align*}
with $\sigma=\int_{-\infty}^{\infty}(\partial_\xi\tilde{\phi}_0)^2 d\xi=\int_{-\infty}^{\infty}(\partial_\xi\tanh(\xi/\sqrt{2}))^2 d\xi=2\sqrt{2}/3$.  Then we integrate the equation \eqref{e:temp0} on $\xi$ in $(-\infty,\infty)$, and notice the matching condition
$$
\lim_{\xi\rightarrow \pm\infty} \tilde p=p^{\pm}, \quad
\lim_{\xi\rightarrow \pm\infty} \partial_{\xi}\tilde{\mathbf v}_1=\mathbf{n}\cdot\nabla \mathbf{v}_0^{\pm}.
$$
We are led to
\begin{equation*}
[-p_0 \mathbf{n}+(\bn\cdot\nabla) \mathbf{v}_0]={\mathsf{B}}\sigma\kappa\mathbf{n}.
\end{equation*}
This is the jump condition for pressure and stress~\citep{magaletti2013sharp}.

Combining the above analysis, in the leading order, we obtain the  standard Navier-Stokes equation for two-phase immiscible flow
\begin{equation}\lbl{e:bulkSharp}
\left\{
\begin{array}{ll}
\mathsf{R_e}\Big(\frac{\partial\mathbf{v}_0^{\pm}}{\partial t}+(\mathbf{v}_0^{\pm}\cdot\nabla)\mathbf{v}_0^{\pm}\Big)=
\mathbf{F}  -\nabla{p}_0^{\pm} +{\Delta}{ \mathbf{v}_0^{\pm}},& \hbox{ in }\Omega^{\pm},\\
\nabla\cdot \mathbf{v}_0^{\pm}=0,& \hbox{ in }\Omega^{\pm},\\
\big[\mathbf{v}_0\big]=0, & \hbox{ on } \Gamma,\\
\big[-p_0 \mathbf{n}+(\bn\cdot\nabla) \mathbf{v}_0\big]={\mathsf{B}}\sigma\kappa\mathbf{n},& \hbox{ on }\Gamma,\\
V_n=\mathbf{v}_0\cdot\mathbf{n} ,& \hbox{ on }\Gamma.
\end{array}
\right.
\end{equation}

\subsection{The boundary conditions}
We now consider the sharp-interface limit of 
the  boundary condition \eqref{e:2.2n}.
We  consider three different choices for the relaxation parameter $\mathrm{V_s}=O(\eps^\beta)$,
$\beta=0,-1,-2$. We show that they correspond
 to different boundary conditions in the sharp-interface limit. 

{\bf Case I. $\mathrm{V_s}=O(1)$.} We first assume that $\mathrm{V_s}$
is a constant independent of the Cahn number $\eps$.

{\it Outer expansion. } 
 Far from the moving contact line,
we can use the same outer expansions as in the bulk. 
{By applying} the expansions~\eqref{e:outexpan} to
\eqref{e:2.2n} {and} considering the leading order terms, 
we easily have the Navier slip boundary condition
\begin{align}\lbl{e:NavierSlip}
&\mathsf{l_s}^{-1} ({v}_{0,\tau}-v_{w})=-  \partial_{n}{v}_{0,\tau},
\quad \mathbf{v}_{0}\cdot\mathbf{n}_S=0,  \qquad \hbox{on } \Gamma_S^{\pm}.
\end{align}
Here ${v}_{0,\tau}=\mathbf{v}_0-(\mathbf{v}_{0}\cdot\mathbf{n}_S)\mathbf{n}_S$ is the tangential velocity,
and $\Gamma_S^{\pm}=\Gamma_{S}\cap\partial\Omega^{\pm}$ is  the boundary of $\Omega^{\pm}$ on
the solid surface $\Gamma_S$.
For the chemical potential $\mu$, we  have
\begin{equation}\label{e:temp_bnd}
\partial_n \mu^{\pm}_0=0.
\end{equation} 
Notice that $\mu_0^{\pm}=0$ and $\phi^{\pm}_0=\pm 1$ in $\Omega^\pm$,
we easily obtain
$$
\phi_0^{\pm}=\pm 1,\quad \mu_0^{\pm}=0, \qquad \hbox{on } \Gamma_S^{\pm}.
$$ 

{\it Inner expansion. } We consider the boundary condition near the contact line. 
Here we denote  the out normal of the solid surface $\Gamma_S$  by $\mathbf{n}_S$,
the  normal of the contact line in tangential surface of $\Gamma_S$  by $\mathbf{m}$ (as shown in Fig.~\ref{fig:1}).
Since the functions at the contact line need to be matched to the outer expansions on $\Gamma_S^{\pm}$
 and also to the expansions inside the domain $\Omega$, it is convenient to
 introduce a different inner expansion near the contact line 
 $L:=\{x\in\Gamma_S\,|\,\phi_\eps=0\}$ as follows. Near the contact line, we introduce two fast {changing} parameters,
\begin{equation*}
\varrho=\frac{(x-x_0)\cdot\mathbf{m}}{\eps}, \qquad \zeta=\frac{(x-x_0)\cdot\mathbf{n_S}}{\eps},
\end{equation*}
with $x_0\in L$. For any function $f(x)$, near $x_0$, it can be written as a function in $(x,\rho,\zeta)$ as  
\begin{equation}
f(x)=\hat{f}(x,\varrho,\zeta).
\end{equation}
The derivative of  $f$ is then rewritten as 
\begin{align*}
\nabla f& =\nabla \hat{f}+\eps^{-1} \mathbf{n}_S\partial_\zeta\hat{f} +\eps^{-1} \mathbf{m}\partial_\varrho\hat{f},\\
\Delta f&=\Delta\hat{f}+\eps^{-1}(\nabla\cdot \mathbf{n}_S\partial_\zeta \hat{f}+\nabla\cdot \mathbf{m}\partial_\varrho \hat{f}
+2\partial_{n\zeta}\hat{f}+2\partial_{m\varrho}\hat{f})+\eps^{-2}(\partial_{\zeta\zeta}\hat{f}+\partial_{\varrho\varrho}\hat{f}),
\end{align*}
and the boundary derivative of $f$ is given by
\begin{align*}
{\partial_n f}=\partial_n \hat{f}+\eps^{-1}\partial_{\zeta} \hat{f},\qquad
\partial_{\tau} f= \partial_{\tau} \hat{f}+\eps^{-1}\partial_{\varrho}\hat{f}\mathbf{m}.
\end{align*}
We also have 
$$\partial_t f=\partial_t \hat{f}-\eps^{-1}\partial_{\varrho}\hat{f}\partial_t x_0\cdot\mathbf{m}. $$
Similarly as before, we  assume that
\begin{align*}
\hat{\mathbf{v}}_{\eps}&=\hat{\mathbf{v}}_0+\eps \hat{\mathbf{v}}_1+\eps^2 \hat{\mathbf v}_2+\cdots,\\
\hat \phi_{\eps}&=\hat \phi_0+\eps \hat \phi_1+\eps^2 \hat \phi_2+\cdots.
\end{align*}
Direct computations give
\begin{equation*}
\hat \mu_\eps=\eps^{-1}\hat \mu_0+\hat \mu_1+\cdots,
\end{equation*}
with
$
\hat\mu_0=-(\partial_{\zeta\zeta}+\partial_{\varrho\varrho})\hat{\phi}_0-\hat \phi_0+\hat \phi_0^3, 
$
and 
$$
\hat{\mathcal{L}}=\hat{\mathcal{L}}_0+\eps \hat{\mathcal{L}}_1+\cdots,
$$
with $\hat{\mathcal{L}}_0=\partial_{\zeta}\hat{\phi}_0+\frac{\partial\gamma_{wf}}{\partial\phi}(\hat{\phi}_0)$.

Substitute the expansions into the CHNS {system} \eqref{e:2.1n} and the boundary condition \eqref{e:2.2n}. In the leading order 
we have the following equation
\begin{equation}\lbl{e:innerbnd}
\left\{
\begin{array}{ll}
(\partial_{\zeta\zeta}+\partial_{\varrho\varrho})\hat{\mu}_0=0, & \hbox{} \\
{\mathsf{B}}\hat{\mu}_0(\partial_{\zeta}\hat{\phi}_0 \mathbf{n}_S+\partial_{\varrho}\hat{\phi}_0\mathbf{m})+(\partial_{\zeta\zeta}+\partial_{\varrho\varrho})\hat{\mathbf{v}}_0=0, & \hbox{}\\
\mathbf{n}_S\cdot\partial_{\zeta}\hat{\mathbf{v}}_0 +\mathbf{m}\cdot\partial_{\varrho}\hat{\mathbf{v}}_0=0, & \hbox{}
\end{array}
\right.
\end{equation}
and the  boundary condition on $\Gamma_S$:
\begin{equation}\lbl{e:innerbndGamma}
\left\{
\begin{array}{ll}
 -\partial_t x_0\cdot\mathbf{m}+ {\hat{v}}_{0,\tau}\cdot\mathbf{m}
=0, & \hbox{} \\
\eps\mathrm{l_s}^{-1} (\hat{v}_{0,\tau} -v_w) +\eps \partial_n \hat{v}_{0,\tau}+  \partial_{\zeta}\hat{v}_{0,\tau}= \mathsf{B}\hat{\mathcal{L}}_0\partial_{\varrho}\hat{\phi}_0, & \hbox{}\\
\partial_\zeta \hat{\mu}_0=0,\qquad \hat{\mathbf{v}}_0\cdot\mathbf{n}_S=0. & \hbox{}
\end{array}
\right.
\end{equation}
Here 
in the second equation of \eqref{e:innerbndGamma}, we keep the terms $\eps\mathrm{l_s}^{-1} (\hat{v}_{0,\tau}-v_w)+\eps \partial_n \hat{v}_{0,\tau}$,
since the slip velocity in the vicinity of the moving contact line might be large compared with the out region~\citep{QianWangSheng2004}. 

{\it Matching condition.} 
We  have the  matching condition
\begin{align*}
&\lim_{\zeta \rightarrow+\infty}\hat{f}=\lim_{x\rightarrow x_0}\tilde f(x,\xi),\\
&\lim_{\varrho \rightarrow\pm\infty}\hat{f}=\lim_{x\rightarrow x_0}f^{\pm}(x).
\end{align*}

In the following, we will use these equations to derive the condition for moving contact lines.
By the matching relation for $\mu$,
\begin{equation*}
\lim_{\zeta \rightarrow+\infty}\hat{\mu}_0=\tilde{\mu}_0=0, \qquad
\lim_{\varrho \rightarrow\pm\infty}\hat{\mu}_0=\mu^{\pm}=0.
\end{equation*}
The first equation of \eqref{e:innerbnd} implies that
\begin{equation}\lbl{e:hatmu0}
\hat{\mu}_0=0.
\end{equation}
This means 
\begin{equation}\lbl{e:phi0CL}
-(\partial_{\zeta\zeta}+\partial_{\varrho\varrho})\hat{\phi}_0-\hat \phi_0+\hat \phi_0^3=0.
\end{equation}
We also have the matching condition for $\hat{\phi}_0$,
\begin{equation*}
\lim_{\zeta \rightarrow+\infty}\hat{\phi}_0=\tilde{\phi}_0(\xi), \qquad
\lim_{\varrho \rightarrow\pm\infty}\hat{\phi}_0=\phi_0^{\pm}=\pm 1.
\end{equation*}
 By \eqref{e:Equprofile}, it is easy to see that
$$
\hat{\phi}_0(\zeta,\varrho)=\tilde{\phi}_0(\xi)
$$
satisfies the equation \eqref{e:phi0CL} and the matching condition
if the relation $$\zeta=\xi/\cos\theta_d,\quad\varrho=\xi/\sin\theta_d$$ holds.
Here $\theta_d$ is the dynamic contact angle which is equal to the angle between 
$\mathbf{n}$ and $\mathbf{m}$ (see Figure~\ref{fig:1}).
This leads to the following relation
\begin{align}\label{e:temp_inner}
\partial_\varrho \hat{f}=\partial_{\xi} \tilde f \sin\theta_d, \quad
\partial_\zeta \hat{f}=\partial_{\xi} \tilde f \cos\theta_d.
\end{align}
In addition, \eqref{e:hatmu0} and the second equation of \eqref{e:innerbndGamma} implies that
$$
(\partial_{\zeta\zeta}+\partial_{\varrho\varrho})\hat{\mathbf{v}}_0=0.
$$
We also have the matching condition
$
\lim_{\zeta \rightarrow+\infty}\hat{\mathbf{v}}_0=\lim_{x\rightarrow x_0}\tilde{\mathbf{v}}_0=\mathbf{v_0}(x_0)
$
and the boundary condition $\hat{\mathbf{v}}_0\cdot\mathbf{n}_S=0$. It is  
easy to see that 
$$ 
\hat{\mathbf{v}}_0 =\mathbf{v_0}(x_0)
$$
with $$\mathbf{v_0}(x_0)\cdot\mathbf{n}_S=0$$ is a solution of the above equation, i.e. $\hat{\mathbf{v}}_0$ is independent of $\zeta$ and $\varrho$.

Using the above relations, we will derive the condition of moving contact line.
The first  equation of \eqref{e:innerbndGamma} implies that
\begin{equation}
\partial_t x_0\cdot\mathbf{m}= {v}_{\tau}(x_0)\cdot\mathbf{m},
\end{equation}
since $\hat{v}_{\tau}=\hat{\mathbf{v}}_0-(\hat{\mathbf{v}}_0\cdot\mathbf{n}_S)\mathbf{n}_S={v}_{\tau}$.
This implies that the normal velocity of the moving contact line in tangential surface is equal to 
the  fluid velocity in this direction.
The second equation~\eqref{e:innerbndGamma} can be reduced to
$$
\eps\mathrm{l_s}^{-1} ( \hat{v}_{0,\tau} -v_w)+\eps \partial_n \hat{v}_{0,\tau}= \mathsf{B}\hat{\mathcal{L}}_0\partial_{\varrho}\hat{\phi}_0.
$$
Integrate this equation with respect to $\varrho$ and we get
\begin{equation}\lbl{e:temp2}
\int_{-\infty}^{+\infty}\eps\mathrm{l_s}^{-1} ( \hat{v}_{0,\tau} - v_w )+  \eps \partial_n \hat{v}_{0,\tau}d \varrho=
\int_{-\infty}^{+\infty}
\mathsf{B}(\partial_{\zeta}\hat{\phi}_0+\frac{\partial\gamma_{wf}}{\partial\phi}(\hat{\phi}_0))\partial_{\varrho}\hat{\phi}_0
d\varrho.
\end{equation}
The left hand side term is 
\begin{align}\lbl{e:inegGNBC}
\int_{-\infty}^{+\infty}(\eps\mathrm{l_s}^{-1} ( \hat{v}_{0,\tau} - v_w )+ \eps \partial_n \hat{v}_{0,\tau})d \varrho=
\int_{interface}(\mathrm{l_s}^{-1} ( \hat{v}_{0,\tau} - v_w )+  \partial_{n}\hat{v}_{0,\tau})d m.
\end{align}
The right hand side term is 
\begin{align}
\int_{-\infty}^{+\infty}
\mathsf{B}(\partial_{\zeta}\hat{\phi}_0+\frac{\partial\gamma_{wf}}{\partial\phi}(\hat{\phi}_0))\partial_{\varrho}\hat{\phi}_0
d\varrho
&=\int_{-\infty}^{+\infty}
\mathsf{B}(\partial_{\xi}\hat{\phi}_0\cos\theta_d+\frac{\partial\gamma_{wf}}{\partial\phi}(\hat{\phi}_0))\partial_{\xi}\hat{\phi}_0d\xi\nonumber\\
&=\mathsf{B}\sigma \cos\theta_d+\mathsf{B}(\gamma_{wf}(1)-\gamma_{wf}(-1))\nonumber\\
&=\mathsf{B}\sigma(\cos\theta_d-\cos\theta_s).
\end{align}
Here we used the Young equation $\sigma\cos\theta_s=\gamma_{wf}(1)-\gamma_{wf}(-1)$
and \eqref{e:temp_inner}. 
Then the equation  \eqref{e:temp2} is reduced to
\begin{equation*}
\int_{interface}(\mathrm{l_s}^{-1} ({v}_{0,\tau} -v_w)+  \partial_{n_S}{v}_{0,\tau})d m=
\mathsf{B}\sigma(\cos\theta_d-\cos\theta_s).
\end{equation*}
This implies that
\begin{equation}
\mathrm{l_s}^{-1} ({v}_{0,\tau} -v_w)+  \partial_{n_S}{v}_{0,\tau}=\mathsf{B}\sigma(\cos\theta_d-\cos\theta_s)\delta_{CL}.
\end{equation}

Combine the above analysis, in the leading order, we have the boundary condition 
\begin{equation}
\left\{
  \begin{array}{ll}
  \mathbf{v}_{0}\cdot\mathbf{n}_S=0, & \\
  V_{CL}={v}_{0,\tau}\cdot\mathbf{m},&\\
  \mathrm{l_s}^{-1} ({v}_{0,\tau} -v_w)+  \partial_{n}{v}_{0,\tau}=\mathsf{B}\sigma(\cos\theta_d-\cos\theta_s)\delta_{CL}. 
  \end{array}
\right.
\lbl{e:sharpbnd}
\end{equation}
This is the sharp-interface version of the generalized Navier slip boundary conditions, which
has been used by~\citet{gerbeau2009,Buscaglia2011Variational,reusken2015finite}.

{\bf Case II.} $\mathsf{V}_s=O(\eps^{-1})$. In  this case, we assume a larger relaxation number 
 $\mathsf{V}_s=\eps^{-1}\alpha $. 
The boundary condition~\eqref{e:2.2n} will be reduced to
\begin{equation}\lbl{e:newbnd}
\left\{
  \begin{array}{ll}
     \frac{\partial {\phi}}{\partial {t}} + {{v}}_{\tau}{\partial_{\tau}}{\phi}=-\eps^{-1}{\alpha} \mathcal{L}(\phi), & \hbox{} \\
  \mathsf{l_s}^{-1} ({v}_{\tau}-v_w)=-  \partial_{n}{v}_{\tau}+ \mathsf{B}\mathcal{L}(\phi)\partial_{\tau}{\phi}, & \hbox{}\\
\nabla\mu\cdot\mathbf{n}_S=0 , \ \ \ \   \mathbf{v}\cdot\mathbf{n}_S=0. & \hbox{}
  \end{array}
\right.
\end{equation}

We repeat the  analysis in Case I to this boundary condition. The only difference is the first equation of \eqref{e:newbnd}.
Using the same inner expansions, the leading order of the first equation of \eqref{e:newbnd} is give by
\begin{align}\lbl{e:innernewbnd0a}
 -\partial_{\varrho}\hat{\phi}_0(\partial_t x_0\cdot\mathbf{m}-\hat{v}_{0,\tau}\cdot\mathbf{m})=\alpha \hat{\mathcal{L}}_0.
\end{align}
We multiply $\partial_{\varrho}\hat{\phi}_0$ to both sides of the equation, and integrate the results 
in $(-\infty,\infty)$. Direct calculations  give 
\begin{equation*}
-(\partial_t x_0\cdot\mathbf{m}-v_{0,\tau}\cdot\mathbf{m})\sin\theta_d=\alpha (\cos\theta_d-\cos\theta_s).
\end{equation*}
This implies 
\begin{equation}\label{e:mcl}
V_{CL}-{v}_{0,\tau}\cdot\mathbf{m}=-\frac{\alpha}{\sin\theta_d} (\cos\theta_d-\cos\theta_s).
\end{equation}

Therefore, the boundary condition in this case is given by
\begin{equation}
\left\{
  \begin{array}{ll}
\mathrm{l_s}^{-1} ( {v}_{0,\tau} - v_w )+  \partial_{n}{v}_{0,\tau}=\mathsf{B}\sigma(\cos\theta_d-\cos\theta_s)\delta_{CL},\quad \mathbf{v}_0\cdot\mathbf{n}=0, & \\
  V_{CL}=v_{0,\tau}\cdot\mathbf{m}-\frac{\alpha}{\sin\theta_d} (\cos\theta_d-\cos\theta_s).&
  \end{array}
\right.
\lbl{e:sharpbnd1}
\end{equation}
Here the velocity of the contact line is different from
 the fluid velocity due to some 
 extra chemical diffusion on the contact line~\citep{Jacqmin2000,Yue2011}.

{\bf Case III.}  $\mathsf{V}_s=O(\eps^{-2})$.  In  this case, we assume 
 $\mathsf{V}_s=\eps^{-2}\alpha$. 
The boundary condition~\eqref{e:2.2n} will be reduced to
\begin{equation} \lbl{e:newbnd1}
\left\{
  \begin{array}{ll}
     \frac{\partial {\phi}}{\partial {t}} + {{v}}_{\tau}{\partial_{\tau}}{\phi}=-\eps^{-2}{\alpha} \mathcal{L}(\phi), & \hbox{} \\
  \mathsf{l_s}^{-1} ({v}_{\tau}-v_w)=-  \partial_{n}{v}_{\tau}+ \mathsf{B}\mathcal{L}(\phi)\partial_{\tau}{\phi}, & \hbox{}\\
\nabla\mu\cdot\mathbf{n}_S=0 , \ \ \ \   \mathbf{v}\cdot\mathbf{n}_S=0. & \hbox{}
  \end{array}
\right.
\end{equation}

Once again, the only difference is the first equation of \eqref{e:newbnd1}. For inner expansions, 
the leading order of the first equation of \eqref{e:newbnd} is give by
\begin{equation}
\alpha\hat{\mathcal{L}}_0(\hat{\phi}_0)=0.
\end{equation}
This leads to 
\begin{equation}
\cos\theta_d=\cos\theta_s.
\end{equation}
The equation implies that the dynamic contact angle is equal to the (static) Young's angle.
Thus, the boundary condition in this case is reduced to
\begin{equation}
\left\{
  \begin{array}{ll}
\mathrm{l_s}^{-1} ( {v}_{0,\tau} - v_w )+  \partial_{n_S}{v}_{0,\tau}=0,\quad \mathbf{v}_0\cdot\mathbf{n}=0, & \\
  \theta_d=\theta_s.&
  \end{array}
\right.
\lbl{e:sharpbnd2}
\end{equation}
The boundary condition is used by~\citet{renardy2001numerical,spelt2005level}. We see that this condition
is correct only for very large wall relaxation case.

\subsection{Summary of the analysis results}
We summarize the main results in this section.
When the mobility parameter $\mathsf{L}_d$ is of order $O(\eps)$, 
the sharp-interface limit of the CHNS system \eqref{e:2.1n} is the standard two-phase flow equation
\begin{equation}\label{e:sharpModel}
\left\{
\begin{array}{ll}
\mathsf{R_e}\Big(\frac{\partial\mathbf{v}}{\partial t}+(\mathbf{v}\cdot\nabla)\mathbf{v}\Big)=
\mathbf{F}  -\nabla{p} +{\Delta}{ \mathbf{v}}& \hbox{ in }\Omega^{\pm},\\
\nabla\cdot \mathbf{v}=0,& \hbox{ in }\Omega^{\pm},\\
\big[\mathbf{v}\big]=0, & \hbox{ on } \Gamma,\\
\big[-p_0 \mathbf{n}+(\bn\cdot\nabla) \mathbf{v}_0\big]={\mathsf{B}}\sigma\kappa\mathbf{n},& \hbox{ on }\Gamma,\\
V_n=\mathbf{v}\cdot\mathbf{n} ,& \hbox{ on }\Gamma.
\end{array}
\right.
\end{equation}
where $V_n$ is the normal velocity of the interface of the two-phase flow.

The different choices of the parameter $\mathrm{V}_s$  lead to 
different sharp interface limits for the GNBC:

{\it Case I.} When $\mathsf{V}_s=O(1)$, the sharp-interface limit
of the boundary condition is 
\begin{equation}
\left\{
  \begin{array}{ll}
\mathrm{l_s}^{-1} ({v}_{\tau}-v_w)+  \partial_{n}{v}_{\tau}=\sigma(\cos\theta_d-\cos\theta_s)\delta_{CL},\quad \mathbf{v}\cdot\mathbf{n}_S=0, &\\
  V_{CL}={v}_{\tau}\cdot\mathbf{m}.&
  \end{array}
\right.
\lbl{e:bnd1}
\end{equation}
The first equation of \eqref{e:bnd1} is the sharp-interface version of the generalized Navier slip boundary condition. It can be understood in the following way~\citep{QianWangSheng2003}:
$$\mathrm{l_s}^{-1} ({v}_{\tau}-v_w)= - \partial_{n}{v}_{\tau}-\frac{1}{\eta}\sigma_Y,$$
where $\sigma_Y$ is the unbalanced Young stress, satisfying
\begin{equation}
-\frac{1}{\eta}\int_{interface}\sigma_Y=\sigma(\cos\theta_d-\cos\theta_s).
\end{equation}
As shown in the MD simulations by~\citet*{QianWangSheng2003}, the unbalance Young force might
leads to near complete slipness of the fluid in the vicinity of the contact line.

The second equation of \eqref{e:bnd1} implies that the velocity of the contact line is equal to the fluid velocity.

{\it Case II.} When $\mathsf{V}_s=O(\eps^{-1})$, the sharp interface limit
of the GNBC is 
\begin{equation}
\left\{
  \begin{array}{ll}
\mathrm{l_s}^{-1} ({v}_{\tau} - v_w)+  \partial_{n}{v}_{\tau}=\sigma(\cos\theta_d-\cos\theta_s)\delta_{CL},\quad \mathbf{v}\cdot\mathbf{n}=0, & \\
  V_{CL}=v_{\tau}\cdot\mathbf{m}-\frac{\alpha}{\sin\theta_d} (\cos\theta_d-\cos\theta_s),&
  \end{array}
\right.
\lbl{e:bnd2}
\end{equation}
where $\alpha = \eps\mathsf{V}_s$ is a constant.
The first equation in \eqref{e:bnd2} is the same as the previous case. 
The second equation in \eqref{e:bnd2} implies that the motion of the contact line is not only
determined by the fluid velocity, but also  by
the chemical diffusion on the boundary.  

From the second equation of \eqref{e:bnd2},  when $\theta_d$ does not change much from $\theta_s$,
we have $$\cos\theta_d-\cos\theta_s\approx-(\sin \theta_d)(\theta_d-\theta_s)+h.o.t.$$
Then the second equation implies that
$ V_{CL}=v_{\tau}\cdot\mathbf{m}+\alpha (\theta_d-\theta_s).$
This implies that $V_{CL}\propto  (\theta_d-\theta_s)$, which is similar to 
the boundary condition derived by~\citet{ren2007boundary}.

{\it Case III.} When $\mathsf{V}_s=O(\eps^{-2})$, the sharp interface limit
of the boundary condition is 
\begin{equation}
\left\{
  \begin{array}{ll}
\mathrm{l_s}^{-1} ({v}_{\tau} - v_w)+  \partial_{n}{v}_{\tau}=0,\quad \mathbf{v}\cdot\mathbf{n}=0, & \\
\theta_d=\theta_s.&
  \end{array}
\right.
\lbl{e:bnd3}
\end{equation}
The boundary condition is  different from the previous two cases. Here the standard 
Navier slip boundary condition is used on the solid boundary and
the dynamic contact angle is  equal to the Young's angle. This boundary condition
has been used by~\citet{spelt2005level}. 

From the above analysis, we have shown the sharp interface limits (in leading order)
for the CHNS equation with the GNBC. For different
choices of the relaxation parameter,
we obtain some different boundary conditions for moving contact lines. 
In applications,
one could choose the parameters according to ones' own purpose. 
{We} would like to remark that 
{we did not consider the effects of different scalings of slip length $l_s$ with respect to $\eps$, and} 
the analysis does not show the convergence rate
of the  sharp-interface limits, which might be important in real applications.
In next section, we will do numerical simulations for the various choices of the
mobility parameter and the relaxation parameter {to} verify the analytical results and {investigate} the convergence 
{rates} in these cases.


\section{Numerical experiments}
We consider a two-dimensional Couette flow 
in a rectangular domain $\Omega = [0, L_x]\times [-1,1]$ with $L_x = 6$. 
The plates on the top and bottom boundaries move in opposite directions with velocity $\bv_w = (\pm 1, 0)$.
We initiate the phase field as
  \begin{eqnarray}\label{eq:initphi}
 \phi (x, y, t=0) = \tanh \Big( \frac{1}{ \sqrt{2}\eps }\big(0.25L_x - | x - 0.5L_x | \big)\Big) . 
  \end{eqnarray}
We set initial velocity  $\bv_0 = (y, 0)$ for the Couette flow.
$\eps$ is gradually reduced to check the convergence of the solution with respect to $\eps$. The values of $Re, B, \mathsf{l_s}$ are fixed as
\begin{eqnarray}\label{para:set0}
    Re = 0.0001, \quad 
    B = 50, \quad
	\mathsf{l_s} = 0.01.
 \end{eqnarray}

We numerically verify the convergence behavior of the diffuse interface model \eqref{e:2.1n}-\eqref{e:2.2n}  for different scalings of $\mathsf{L_d}, \mathsf{V_s}$  with respect to $\eps$, using a second
order time marching scheme coupled with a spectral method for spatial variables recently proposed by~\citet{YangYu2017}. For clarity, we also list the algorithm in the appendix. 

Fig.~\ref{fig:T02} shows the snapshots of the two-phase interface at $T=0.2$ in the simulation results of Couette flow using different  mobility parameter $\mathsf{L_d}$ and relaxation parameter $\mathsf{ V_s}$. In these experiments, we set $\theta_s=90^\circ$.
Due to symmetry, we only show the bottom part of the left interface in each case.

\begin{figure}
 \centering
 	\subfigure[]{
	\includegraphics[width=0.28\textwidth,height=0.18\textheight]{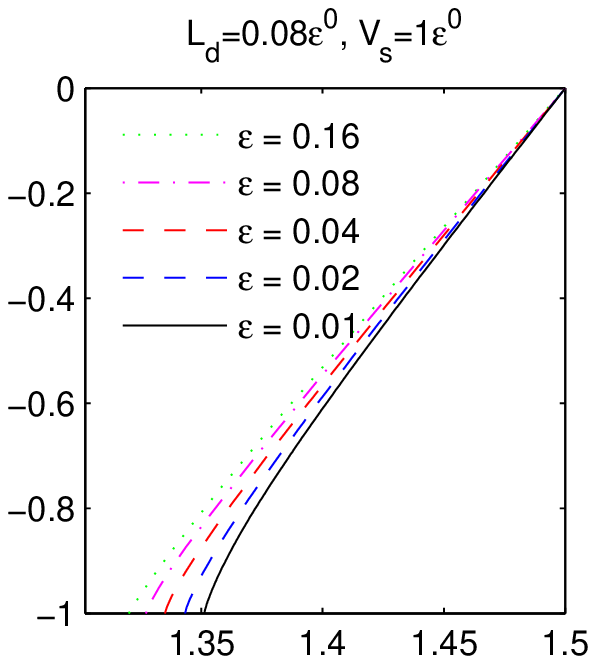}
	}
	\subfigure[]{
		\includegraphics[width=0.28\textwidth,height=0.18\textheight]{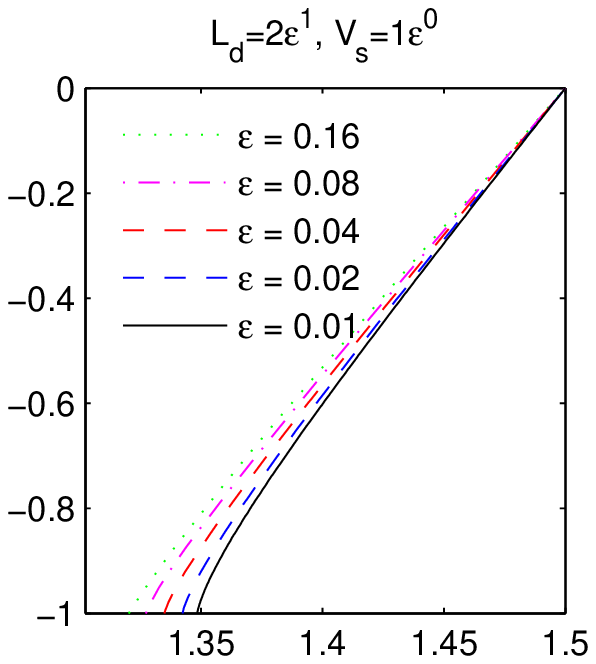}
	}
	\subfigure[]{
		\includegraphics[width=0.28\textwidth,height=0.18\textheight]{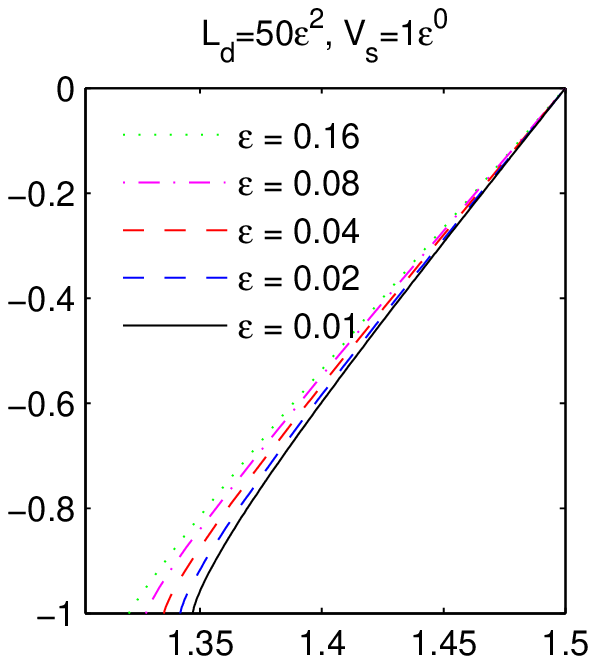}
	}\\
		
     \subfigure[]{
	\includegraphics[width=0.28\textwidth,height=0.18\textheight]{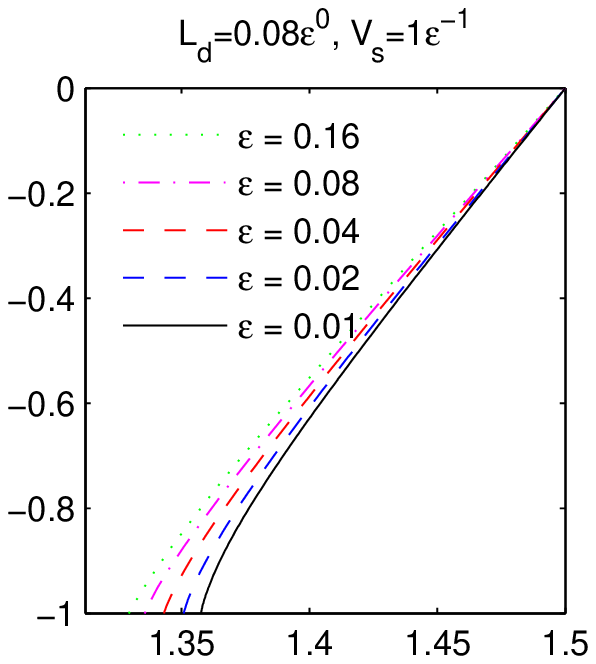}
	}
    \subfigure[]{
	 	\includegraphics[width=0.28\textwidth,height=0.18\textheight]{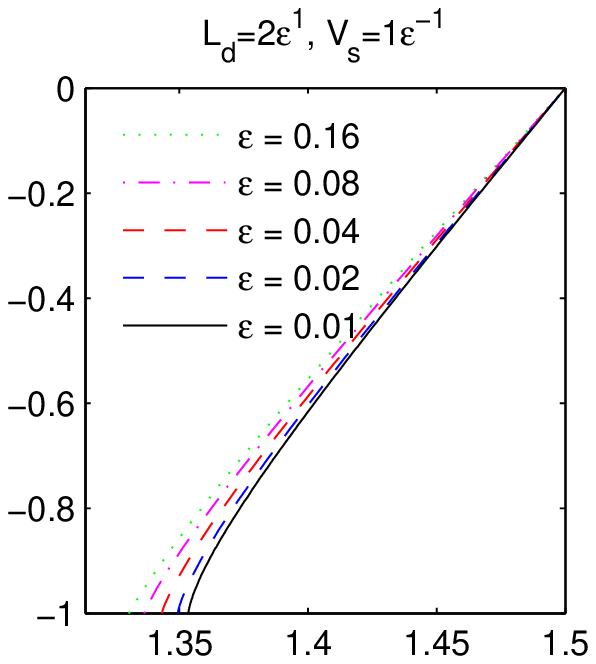}
	}
    \subfigure[]{
	     \includegraphics[width=0.28\textwidth,height=0.18\textheight]{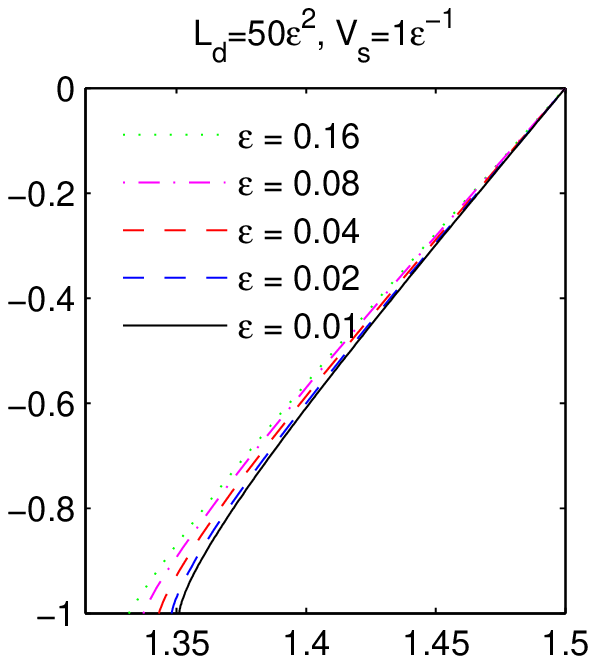}
	 }\\

    \subfigure[]{
	\includegraphics[width=0.28\textwidth,height=0.18\textheight]{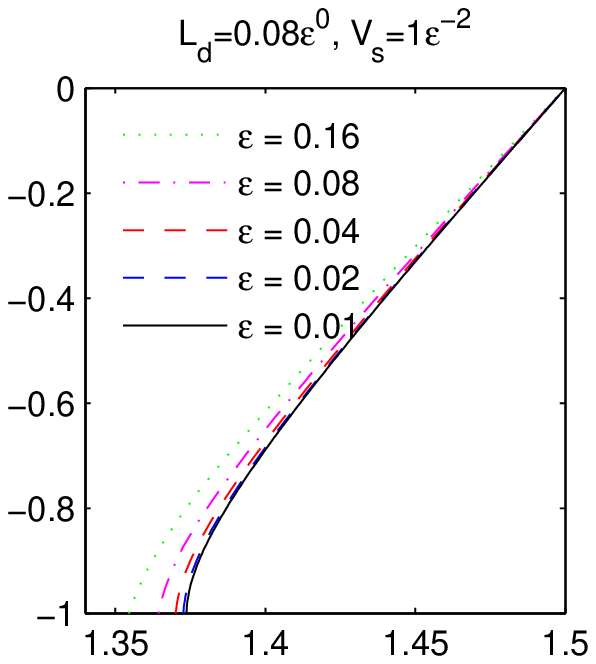}
	}
	 \subfigure[]{
	 	\includegraphics[width=0.28\textwidth,height=0.18\textheight]{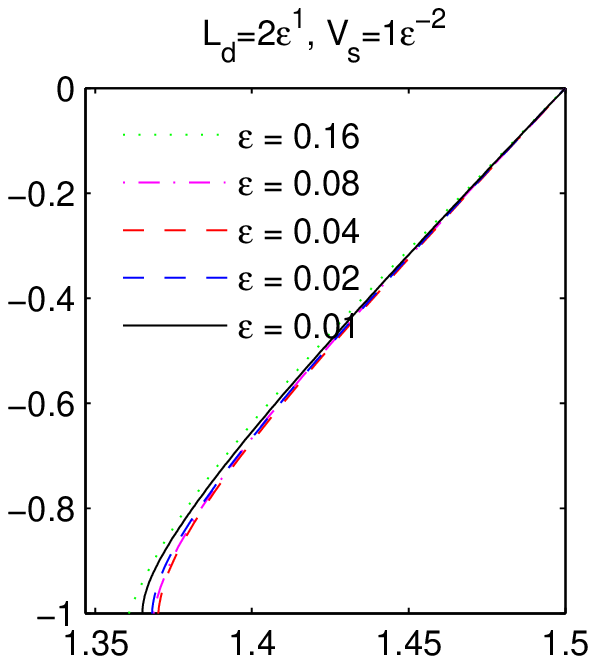}
	 }
	 \subfigure[]{
	 	\includegraphics[width=0.28\textwidth,height=0.18\textheight]{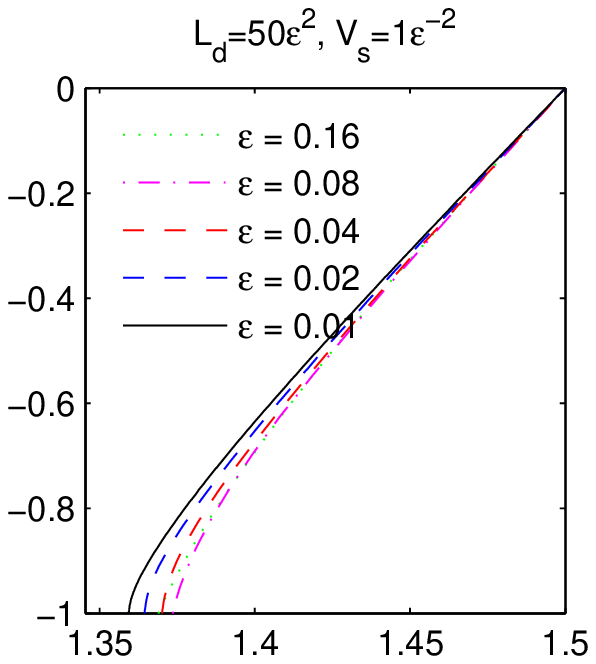}
	 }\\

    \subfigure[]{
	\includegraphics[width=0.28\textwidth,height=0.18\textheight]{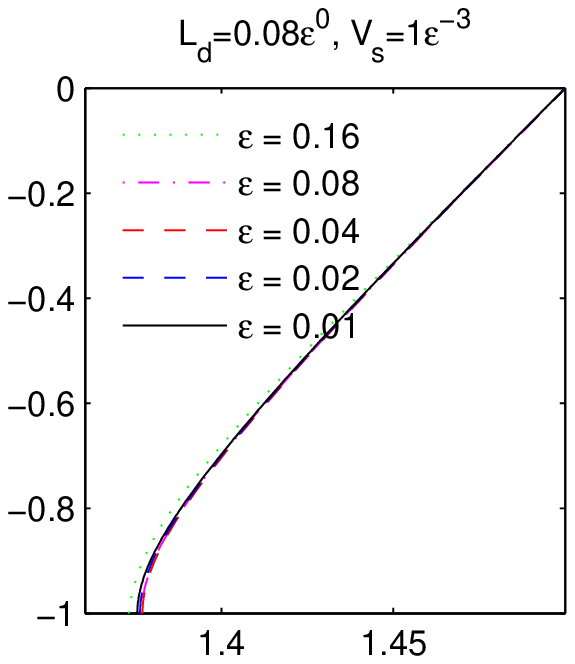}
	}
	 \subfigure[]{
	\includegraphics[width=0.28\textwidth,height=0.18\textheight]{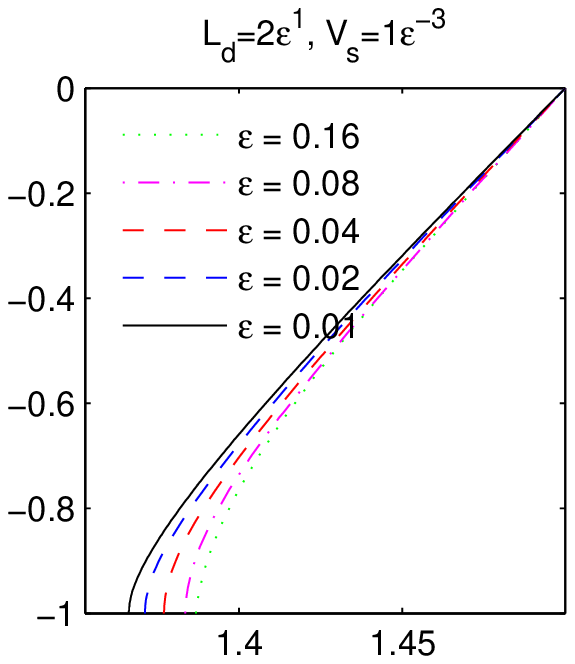}
	}
		 \subfigure[]{
	\includegraphics[width=0.28\textwidth,height=0.18\textheight]{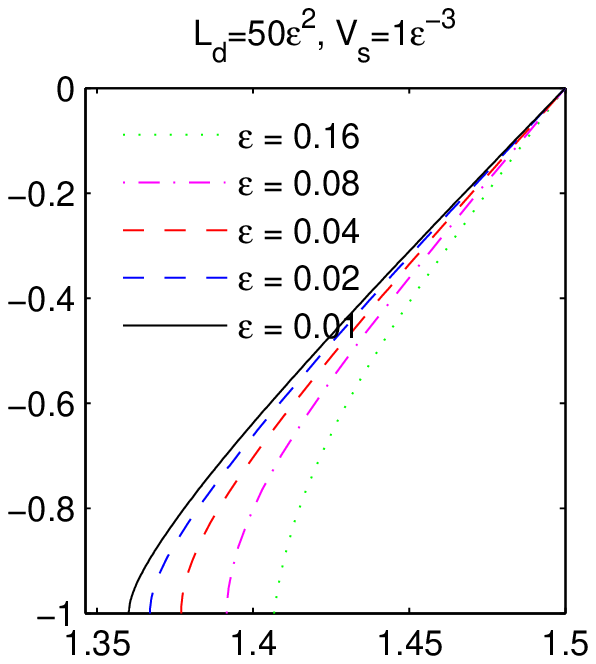}
	}
 	\caption{\label{fig:T02}Numerical results at $T=0.2$ with different $\mathsf{L_d}$ and $\mathsf{V_s}$ values, $\theta_s=90^\circ$. Since the contact lines are symmetric with respect to point $(1.5,0)$, we only plot the bottom parts to show the convergence.
	} 
 \end{figure}
The {left column} of Fig.~\ref{fig:T02} show the results for $\mathsf{L_d}=O(1)$ and $\mathsf{V_s}=O(\eps^\beta)$, with $\beta=0,-1,-2,-3$, respectively. 
It is known that for this case, the Navier-Stokes-Cahn-Hilliard system converges to 
coupled Navier-Stokes  and Hele-Shaw equations~\citep{wang2007sharp}. 
We could also see that the dynamic contact angle approaches to the Young's angle with increasing
relaxation parameter $\mathsf{V_s}$. 
In the largest relaxation 
parameter $\mathsf{V_s}=O(\eps^{-3})$ case, the dynamic contact angle is almost equal to the Young's angle.
We observe that this case {exhibits} the best convergence rate to a sharp-interface limit.
The results are consistent with the numerical observations by~\citet*{Yue2010}.
 In their experiments,  the boundary condition 
$\mathcal{L}(\phi)=0$ {is used}, which corresponds to a infinite large parameter $\mathsf{V}_s$. 
They found that  the sharp interface limit is obtained only when the mobility
parameter is of order $O(1)$.

The middle column of Fig.~\ref{fig:T02} show the results for $\mathsf{L_d}=O(\eps)$ and $\mathsf{V_s}=O(\eps^\beta)$, with $\beta=0,-1,-2,-3$, respectively. 
For different choice of $\mathsf{V_s}$, the sharp-interface limits of the diffuse interface model are slightly different.
With increasing relaxation parameter $\mathsf{V_s}$ (or decreasing $\beta$), the dynamic contact angle will approach to the stationary contact angle $\theta_s=90^\circ$. From Fig.~\ref{fig:T02} (h) and (k), we could see that the dynamic contact angle is almost $90^\circ$ for small $\eps$. 
These results are consistent with the asymptotic analysis in the previous section. 
More interestingly, the different choices of $\mathsf{V_s}$ might also affect 
the convergence rates. In seems that the convergence rate to the sharp-interface
limit for $\mathsf{V_s}=O(\eps^{-1})$ is slightly better than  other cases.
For the case $\mathsf{V_s}=O(\eps^{-3})$, the convergence rate seems  very slow. 
This  indicates the phase-field model might not convergence for infinite large relaxation parameter when $\mathsf{L_d}=O(\eps)$,
as shown by~\citet*{Yue2010}.

 The  {right column} of Fig.~\ref{fig:T02} show the results for $\mathsf{L_d}=O(\eps^{2})$ and $\mathsf{V_s}=O(\eps^\beta)$, with $\beta=0,-1,-2,-3$, respectively.  In this case, the Navier-Stokes-Cahn-Hilliard system  still converges to that
 standard incompressible two-phase Navier-Stokes equations. 
 ~\citet{magaletti2013sharp} considered the case without moving contact lines, found that $\mathsf{L_d}=O(\eps^{2})$ give the best convergence rate.
 Here we show some numerical results for moving contact line problems.  
In this case, the choice of $\mathsf{V_s}=O(\eps^{-1})$ seems correspond to
slightly faster convergence rate than  other choices. {This is similar to the $\mathsf{L_d}=O(\eps)$ case.} 
 
{On the other hand}, we  observe that when the boundary relaxation $\mathsf{V_s}=O(\eps^{\beta})$ with $\beta\geq -1$ {(the first and second row of Fig.~\ref{fig:T02})},
the boundary diffusion can be considered as small, then the scaling $L_d=O(\eps^2)$ 
gives the best convergence rate. This is similar to the results of
phase-field model without contact lines obtained by~\citet{magaletti2013sharp},
where they give an elaborate analysis.
On the other hand, when $V_s=O(\eps^{\beta})$ with $\beta\leq -2$ {(the third and fourth row of Fig.~\ref{fig:T02})}, the boundary diffusion
is relatively large, then $L_d=O(1)$ 
gives best convergence rate. This is consistent to the finding by~\citet*{Yue2010}. 
 The observation is helpful in the real applications using the GNBC model.

\begin{figure}
 \centering
 \subfigure[]{
	\includegraphics[width=0.26\textwidth,height=0.21\textheight]{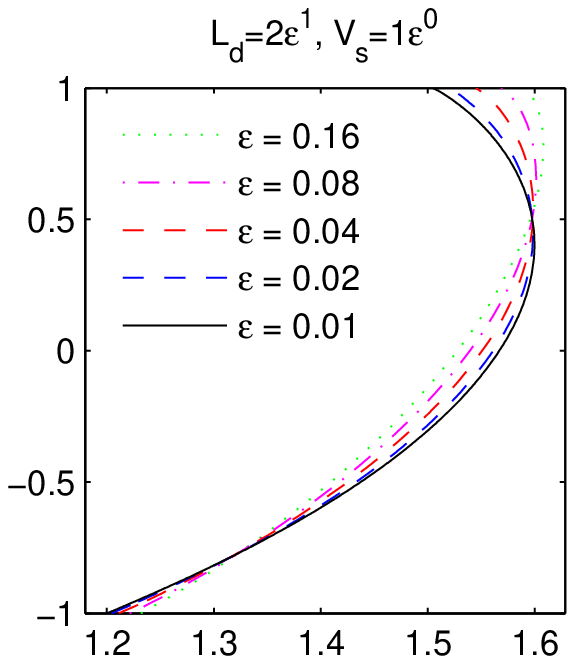}
	}
  \subfigure[]{
  	\includegraphics[width=0.25\textwidth,height=0.21\textheight]{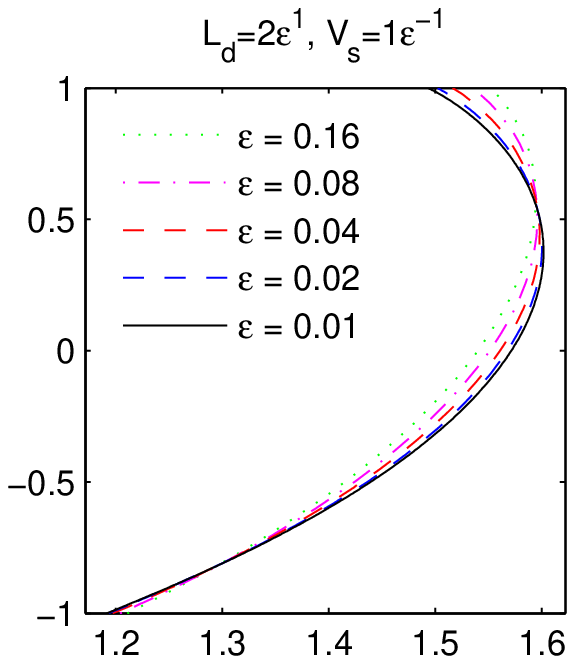}
  }
  \subfigure[]{
 	\includegraphics[width=0.25\textwidth,height=0.21\textheight]{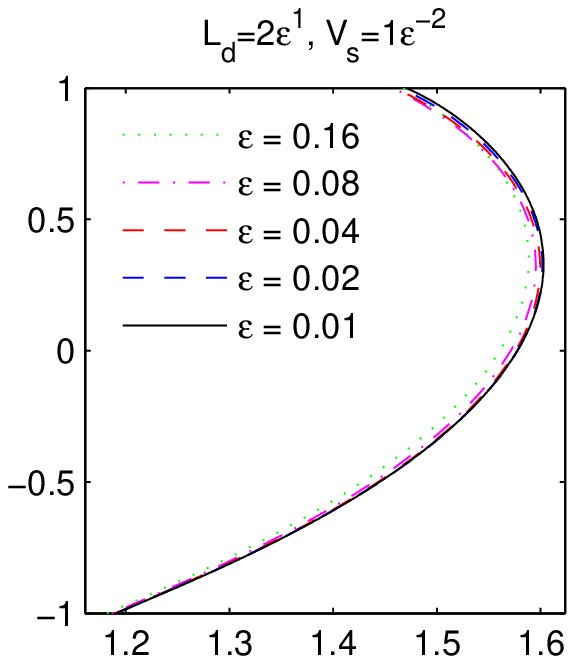}
 }\\

 	\subfigure[]{
	\includegraphics[width=0.26\textwidth,height=0.21\textheight]{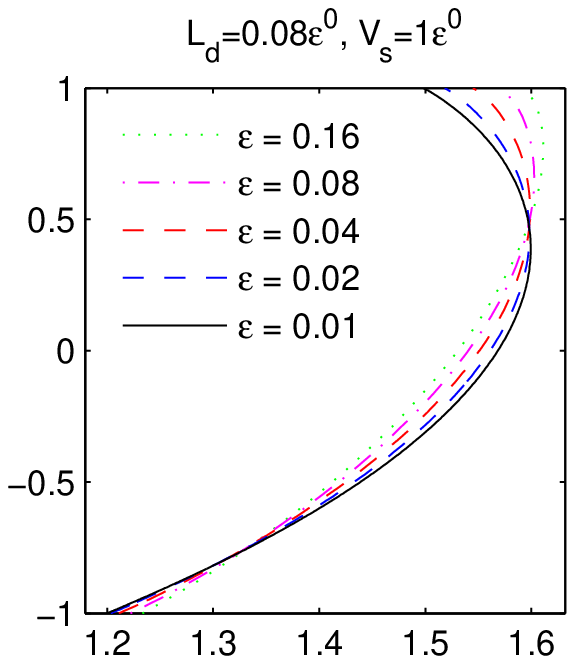}
	} 
	 \subfigure[]{
	\includegraphics[width=0.25\textwidth,height=0.21\textheight]{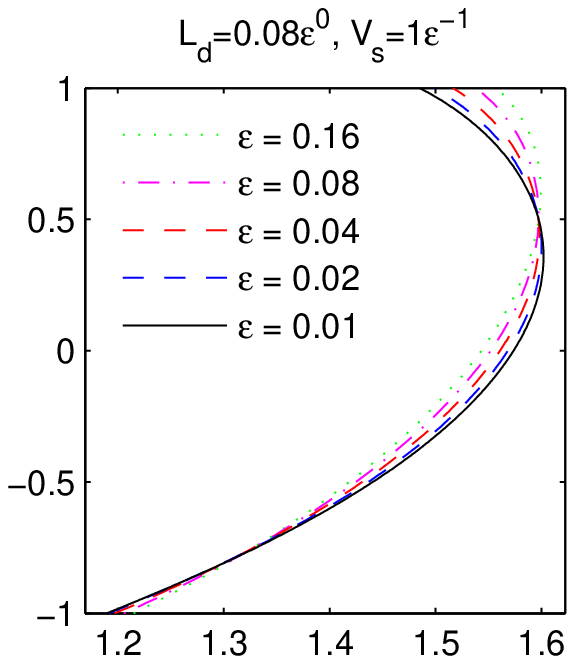}
	}
	 \subfigure[]{
	 	\includegraphics[width=0.25\textwidth,height=0.21\textheight]{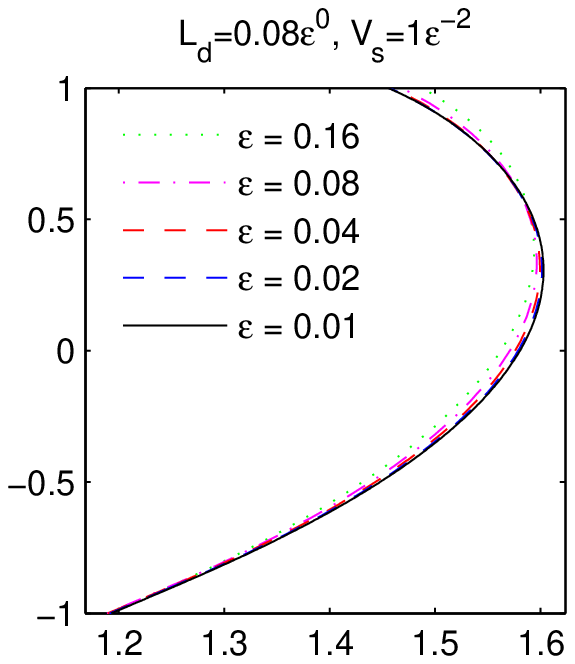}
	 } 	
	 \caption{\label{fig:T02Vs12} Numerical results at $T=0.2$  with different $\mathsf{L_d}$ and $\mathsf{V_s}$ values, $\theta_s=60^\circ$.} 
\end{figure}  
We also did experiments for the case when $\theta_s=60^\circ$. 
The numerical results are similar to the case when $\theta_s=90^\circ$.
In Fig.~\ref{fig:T02Vs12} we present only a few snapshots of MCLs at $T=0.2$ for the case when $\theta_s=60^\circ$.
Here we focus on the differences between the choices $\mathsf{L_d}=O(\eps)$ and $\mathsf{L_d}=O(1)$.
  We show three cases $\mathsf{V_s}=O(1)$, $O(\eps^{-1})$ and  $O(\eps^{-2})$.
From the figure, we could see that the convergence rate for  $\mathsf{L_d}=O(\eps)$ is better than $\mathsf{L_d}=O(1)$ when
$\mathsf{V_s}$ is $O(1)$ and $O(\eps^{-1})$.

\section{Conclusion}

We studied the convergence behavior with respect to the Cahn number $\eps$ of a phase field moving contact line model incorporating dynamic contact line condition in different situations, using asymptotic analysis and numerical simulations. In particular, we considered the situations of 
$\eps$-dependent mobility $\mathsf{L_d}$ and boundary relaxation $\mathsf{V_s}$. This extends the study by~\citet*{Yue2010} and~\citet*{wang2007sharp}. 
\citet*{Yue2010} showed that $\mathsf{L_d}=O(1)$ is the only proper choice for the sharp-interface limit of a diffuse interface model with no slip boundary condition. 
~\citet{wang2007sharp} deduced that the sharp-interface limit of the phase-field model with the GNBC for the case $\mathsf{L_d}=O(1)$ obeys a Hele-Shaw flow. 

We did asymptotic analysis for the phase-field model with the GNBC for the case $\mathsf{L_d}=O(\eps)$. We show that the sharp-interface limit is the  incompressible two-phase Navier-Stokes equations
with standard jump condition for velocity and stress on the interfaces. We also 
show that the different choices of the scaling of $\mathsf{V_s}$ correspond 
to different boundary conditions in the sharp-interface limit.

Our numerical results show that when the boundary relaxation $\mathsf{V_s}=O(\eps^{\beta})$ with $\beta\geq -1$,
the boundary diffusion can be considered as small, the scaling $L_d=O(\eps^\alpha), \alpha=0,1,2$ all give 
convergence, but $\alpha=2$ gives the best convergence rate. This is consistent to the results of
phase-field model without contact lines obtained by~\citet{magaletti2013sharp}.
On the other hand, when $V_s=O(\eps^{\beta})$ with $\beta\leq -2$, $L(d)=O(\eps^\alpha), \alpha=0,1$ will give better convergence rate, while $\alpha=2$ also exhibits convergence. The case $\alpha=0$ 
give best convergence rate for $\mathsf{V_s}$ very large is consistent to the finding by~\citet*{Yue2010}. 
The larger convergence region of $\alpha$ is due to the fact that the generalized Navier slip boundary condition and the dynamic contact line condition are incorporated in the phase field MCL model we used.
Our analysis and numerical studies will be  helpful in the real applications using the GNBC model.

\section*{Appendix: The numerical scheme}

The CHNS system \eqref{e:2.1n}-\eqref{e:2.2n} are solved using a second-order  accurate and energy stable time marching scheme basing on invariant energy quadratization skill developed by~\citet{YangYu2017}. For clarity, we list the details of the scheme below.
Let $\delta t >0$ be the time step-size and set $t^n=n \delta
t$. For any function $S(\bx, t)$, 
let $S^n$ denotes the
numerical approximation to $S(\cdot,t)|_{t=t^n}$, and
$S_\star^{n+1}:=2S^{n}-S^{n-1}$.
 We introduce $U^0 = (\phi^0)^2-1$, $W^0 = \sqrt{\hat\gamma_{wf}(\phi^0)}$, where
 \begin{equation}
 \hat\gamma_{wf}(\phi) = \left\{ 
 	 \begin{array}{ll}
	 \frac{\sqrt{2}}{3}-\frac{\sqrt{2}}{6}\cos\theta_s(3\phi-\phi^3),& 
	 \mbox{if}\ |\phi|\le 1,\\
	 \frac{\sqrt{2}}{3}-\frac{\sqrt{2}}{3}\cos\theta_s,& \mbox{otherwise}.
 	 \end{array}
 	\right.
 \end{equation}

	Assuming that $(\phi, \bv, p, U, W)^{n-1}$ and $(\phi,
	\bv, p, U, W)^{n}$ are already known, we compute
	$\phi^{n+1}, \bv^{n+1}$, $p^{n+1}, U^{n+1}, W^{n+1}$ in two steps:
	
	{\bf Step 1:} We update $\phi^{n+1}, \tilde\bv^{n+1}, U^{n+1}, W^{n+1}$ as follows,
	\begin{eqnarray}
	&&\frac{3\phi^{n+1}-4\phi^n+\phi^{n-1}}{2\delta t}+
	\Grad\cdot(\tilde\bv^\none\phi_\star^{n+1})=
	\mathsf{L_d}\Delta\mu^\none,\label{bdf:1}\\
	&&\mu^\none=-\veps\Delta\phi^\none+\frac{1}{\veps}\phi_\star^{n+1} U^\none,\label{bdf:2}\\
	&&3U^{n+1}-4U^n+U^{n-1}=2\phi_\star^{n+1}
	(3\phi^{n+1}-4\phi^n+\phi^{n-1}),\label{bdf:3}\\
	&&\mathsf{Re}\Big[\frac{3\wtilde\bv^{n+1}\!\!\!-\!\!4\bv^n\!+\!\bv^{n\!-\!1}}{2\delta t}
	+B(\bv_\star^{n+1}\!,\wtilde\bv^\none)\! \Big]-\!\Delta\wtilde \bv^\none\!\!+\Grad p^n\!\!+\mathsf{B}\phi_\star^{n+1}\Grad\mu^\none=F^{n+1}, \label{bdf:4}
	\end{eqnarray}
	with the boundary conditions
	\begin{eqnarray}
	&&\wtilde\bv^{n+1}\cdot\Bn=0 ,\label{bdfbd:1}\\
	&&\partial_{\Bn}\wtilde\bv_\tau^\none=
	-\mathsf{l_s}^{-1}(\wtilde\bv^\none -\bv_w)-\frac{1}{\mathsf{V_s}}\dot{\phi}^\none \Grad_\tau\phi_\star^{n+1} ,\label{bdfbd:2}\\
	&&\partial_\Bn\mu^\none=0 ,\label{bdfbd:3}
	\end{eqnarray}
	\begin{eqnarray}
	&&\veps\partial_\Bn\phi^\none=-\frac{1}{\mathsf{V_s}}\dot{\phi}^\none-Z(\phi_\star^{n+1}) W^\none,\label{bdfbd:4}\\
	&&3W^{n+1}-4W^n+W^{n-1}=\frac 12 Z(\phi_\star^{n+1})(3\phi^{n+1}-4\phi^n+\phi^{n-1}), \label{bdfbd:5}
	\end{eqnarray}
	where 
	\begin{eqnarray}
&B(\bu, \bv)=(\bu\cdot\Grad)\bv+\frac 12(\Grad\cdot\bu)\bv,\\
&\dot{\phi}^\none=\frac{3\phi^{n+1}-4\phi^n+\phi^{n-1}}{2\delta t}+\wtilde\bv_\tau^\none\cdot\Grad_\tau
	\phi_\star^{n+1},\\
&Z(\phi)={\hat\gamma'_{wf}(\phi)}/\sqrt{\hat\gamma_{wf}(\phi)}.
	\end{eqnarray}

	{\bf Step 2:} We update $\bv^{n+1}$ and $p^{n+1}$ as follows, 
	\begin{eqnarray}
	&&\frac{3\mathsf{Re}}{2\delta t}\big(\bv^{n+1}-\wtilde\bv^{n+1}\big)
	+\Grad(p^{n+1}-p^n)=0,\label{bdf:5}\\
	&&\Grad\cdot\bv^{n+1}=0,\label{bdf:6}
	\end{eqnarray}
	with the boundary condition
	\begin{eqnarray}\label{bdfbd:6}
	\bv^{n+1}\cdot\Bn=0 \quad \mbox{on}\ \Gamma.
	\end{eqnarray}
The above scheme are further discretized in space using an efficient Fourier-Legendre spectral method, see~\citep{shen2015efficient} and~\citep{yu_numerical_2017} for more details about the 
spatial discretization and solution procedure.	
The scheme \eqref{bdf:1}-\eqref{bdfbd:6} can be proved to be unconditional energy stable~\citep{YangYu2017}. But to get 
accurate numerical results, we have to take time step-size small enough. 
In all the simulations in this paper we use $\delta t= 0.0001$ and the first order scheme 
proposed by~\citet*{shen2015efficient} is used to generate the numerical solution at $t=\delta t$ to start up
the second order scheme. 

\section*{Acknowledgments}
We thank Professor Xiaoping Wang and Professor Pingbing Ming for helpful discussions.
The work of H. Yu was partially supported by
NSFC under grant 91530322, 11371358 and 11771439.
The work of Y. Xu was partially supported by NSFC projects 11571354 and 91630208.
The work of Y. Di was partially supported by
NSFC under grant 91630208 and 11771437.

\bibliographystyle{jfm}
\bibliography{literature}

\end{document}